\documentclass[10pt]{article}
\usepackage[nohead,margin=1.0in]{geometry}
\usepackage{amssymb, amsmath, amsthm,amsfonts}
\usepackage{amsmath}
\usepackage{graphicx,epsfig}
\usepackage[tight]{subfigure}
\usepackage{cite}
\usepackage{times}
\usepackage{float}
\graphicspath{{./pics/}}
\usepackage{color}
\usepackage{epstopdf}
\usepackage{mathrsfs}
\usepackage{threeparttable}
\usepackage{algorithm,algorithmic}
\usepackage{booktabs}
\usepackage{verbatim}
\usepackage{bm}
\numberwithin{equation}{section}
\newtheorem{theorem}{Theorem}[section]
\newtheorem{remark}{Remark}[section]

\newtheorem{lemma}{Lemma}[section]

\newtheorem{assumption}[theorem]{Assumption}
\newtheorem{example}{Example}[section]
\allowdisplaybreaks

\def\II{(\Omega)}
\def\Dal{\partial_t^\alpha}

\renewcommand{\d}{{\rm d}}

\title{Regularity Analysis and High-Order Time Stepping Scheme for Quasilinear Subdiffusion\thanks{The work of BJ is supported by Hong Kong RGC General Research Fund (Project 14306423), and a
start-up fund from The Chinese University of Hong Kong. The work of ZZ is supported by Hong Kong
Research Grants Council (15303122) and an internal grant of Hong Kong Polytechnic University (Project
ID: P0038888, Work Programme: 1-ZVX3).}}

\author{Bangti Jin\thanks{Department of Mathematics, The Chinese University of Hong Kong, Shatin, New Territories, Hong Kong, P.R. China (\texttt{bangti.jin@gmail.com, b.jin@cuhk.edu.hk, qmquan@cuhk.edu.hk}).} \and Qimeng Quan\footnotemark[2] \and Barbara Wohlmuth\thanks{Department of Mathematics, SoCIT,
Technical University of Munich, Boltzmannstr. 3,
Garching, 85748, Germany (\texttt{wohlmuth@cit.tum.de})} \and Zhi Zhou\thanks{Department of Applied Mathematics,The Hong Kong Polytechnic University, Kowloon, Hong Kong, P.R. China (\texttt{zhizhou@polyu.edu.hk})}}

\begin{document}

\maketitle

\begin{abstract}
In this work, we investigate a quasilinear subdiffusion model which involves a fractional derivative of order $\alpha \in (0,1)$ in time and a nonlinear diffusion coefficient. First, using smoothing properties of solution operators for linear subdiffusion and a perturbation argument, we prove several pointwise-in-time regularity estimates that are useful for numerical analysis. Then we develop a high-order time stepping scheme for solving quasilinear subdiffusion, based on convolution quadrature generated by second-order backward differentiation formula with correction at the first step. Further, we establish that the convergence order of the scheme is $O(\tau^{1+\alpha-\epsilon})$ without imposing any additional assumption on the regularity of the solution. The analysis relies on refined Sobolev regularity of the nonlinear perturbation remainder and smoothing properties of discrete solution operators. Several numerical experiments in two space dimensions show the sharpness of the error estimate.

\vskip5pt
\noindent\textbf{Keywords}: quasilinear subdiffusion,  pointwise-in-time regularity, convolution quadrature, error estimate
\end{abstract}

\section{Introduction}
Let $\Omega\subset \mathbb{R}^d$ ($d=1,2,3$) be an open bounded domain with a $C^2$ boundary $\partial\Omega$. Consider the following initial-boundary value problem of quasilinear subdiffusion
\begin{equation}\label{eqn:fde}
	\left\{\begin{aligned}
		 \Dal u - \nabla\cdot(a(u)\nabla u) &= f, \ &\mbox{in}&\ \Omega\times(0,T), \\
		u&=0, \ &\mbox{on}&\ \partial\Omega\times(0,T), \\
		u(0)&=u_0, \ &\mbox{in}&\ \Omega,
	\end{aligned}
	\right.
\end{equation}
where $T>0$ is the terminal time and the notation $ \Dal  u$ denotes the Caputo fractional derivative of order $\alpha\in(0,1)$ in time $t$, defined by \cite[p. 92]{kilbas2006theory}:
\begin{equation*}
	 \Dal  u(t):=
		\frac{1}{\Gamma(1-\alpha)}\int_0^t(t-s)^{-\alpha}u'(s)\ {\rm d}s,
\end{equation*}
with $\Gamma(z)=\int_{0}^{\infty}e^{-s}s^{z-1}\ \mathrm{d}s$, $\Re(z)>0$, being Euler's Gamma function. The functions $f$ and $u_0$ are given source and initial data, respectively. The diffusion coefficient $a(v)$ belongs to the admissible set
\begin{equation*}
	\mathcal{A}:=\{C^1(\mathbb{R}) :  a_0\leq a(v) \leq a_1,\ \mbox{with}\ |a'(v)|\leq B \ \mbox{for any}\ v\in\mathbb{R}\},
\end{equation*}
with $0<a_0\leq a_1<\infty$ being the lower and upper bounds on the diffusivity $a(u)$, respectively. When $\alpha=1$, the model recovers the standard quasilinear parabolic problem, whose numerical analysis has been extensively studied; see \cite[Chapter 13]{Thome2006GalerkinFE} and \cite{AkrivisLubich:2015,AkrivisLiLubich:2017} and the references therein.

The subdiffusion model has attracted much attention in recent years due to its remarkable capability to describe anomalously slow diffusion processes, which have been observed in diverse scientific areas in engineering, physics and biology etc. These processes can be accurately described by continuous time random walks (with the waiting time between consecutive jumps following a heavy-tail distribution) at a microscopical level, and the continuous subdiffusion model describes the evolution of the probability density function (in $\mathbb{R}^d$) of particles appearing at time $t$ and spatial location $x$ macroscopically. In subdiffusion processes, the mean squared particle displacement exhibits a sublinear growth with time $t$ (hence the name subdiffusion), deviating from linear growth of normal diffusion.  This class of mathematical models has found many successful applications, e.g., solute transport in heterogeneous media \cite{dentz2004time,berkowitz2006modeling} and protein transport within membranes \cite{Kou:2008,kou2004generalized,ritchie2005detection}.

Over the past two decades, there have been significant progress in developing time discretization schemes and providing numerical analysis of linear subdiffusion; See the recent monograph \cite{JinZhou:2023book} for an overview. Existing time-stepping schemes for approximating the time-fractional derivative can roughly be divided into two groups: piecewise polynomial interpolation (see, e.g.,  \cite{lin2007finite,alikhanov2015new}) and convolution quadrature \cite{Lubich:1986,Lubich:1988}. In the context of polynomial interpolation, the L1 scheme \cite{lin2007finite}, based on linear interpolation, is one of the most popular methods. Under strong solution regularity conditions, it can achieve a local truncation error $O(\tau^{2-\alpha})$. However, in practice, the L1 scheme, as well as other similar schemes, typically exhibits only first-order convergence on uniform grids, due to a lack of the solution regularity. Indeed, subdiffusion typically involves weak solution singularities at $t=0$, even when the initial condition $u_0$ and the source $f$ are smooth \cite{SakamotoYamamoto:2011}. To tackle the challenge, graded or non-uniform meshes have been proposed \cite{Stynes:2017}. These schemes can compensate the initial singularity but the flexibility of the meshes incurs delicate theoretical analysis of high-order schemes; see, e.g., \cite[Sect. 3]{alikhanov2015new} for the L2-1$_\sigma$ formula. Meanwhile, in the context of convolution quadrature (CQ),
one constructs a discrete kernel on a uniform grid and approximates the relevant continuous kernel in the Laplace domain. {The construction preserves the stability properties of the underlying ODE solvers, which greatly facilitates the theoretical analysis of the schemes \cite{Lubich:1986,Lubich:1988} and enables fast and oblivious implementations \cite{SchadleLopezLubich:2006,LubichSchadle:2002,BanjaiLopez:2019,Fischer:2019}.} Jin et al \cite{jin2016two} proposed an initial correction method for the second-order backward differentiation formula (BDF2) for both subdiffusion  and diffusion wave ($1<\alpha<2$), and derived error bounds for the scheme; see  also \cite{jin2017correction} for high-order corrected BDF schemes. The key is to determine corrections at the starting step in order to ensure an $O(\tau^2)$ accuracy. See also \cite{BaffetHesthaven:2017,KhristenkoWohlmuth:2023} for schemes based on rational approximations. {The adaptation of these schemes to the more involved semilinear subdiffusion case has been studied in \cite{jin2018numerical,Kopteva:2020,wangzhou2020high,LiMa:SINUM2022}, with applications in nonlocal-in-time gradient flow models
\cite{FritzKhristenkoWohlmuth:2023, DuYangZhou:2020, LiSalgado:2023, Karaa:2021, TangYuZhou:2019}.}

So far there are few theoretical and numerical studies on quasilinear subdiffusion. This is largely attributed to the following two facts. First, the time-dependent nature and nonlinearity of the diffusion coefficient $a(u)$ render conventional analysis techniques, e.g., Laplace transform and separation of variables, not directly applicable. Second, the strong nonlinearity of the quasilinear elliptic operator poses big challenges to analyzing the solution regularity, even when the diffusion coefficient $a$ is smooth. Zacher \cite{zacher2012global} proved the existence and uniqueness of solutions to problem \eqref{eqn:fde} and established an $L^p$ regularity. Giga and Namba \cite{GigaNamba:2017}  proved the existence of viscosity solutions for a time-fractional Hamilton-Jacobi equations by adapting Perron's method, and also obtained two stability results. Liu  et al \cite{LiuRockner:2018} studied a general class of time-fractional models using the monotonicity method, including quasi-linear subdiffusion.
Numerically, L\'{o}pez-Fern\'{a}ndez  and Plociniczak \cite{lopez2023convolution} proposed a space semidiscrete scheme and a fully discrete semi-implicit scheme using backward Euler convolution quadrature, and derived error estimates for the schemes.
The error analysis in \cite{lopez2023convolution} is based on a new Gronwall's inequality, a novel coercivity estimate for convolution quadrature that enables the use of the energy method, and most importantly, certain assumptions on several \textit{a priori} pointwise-in-time regularity estimates of the solution.
See also \cite{Plociniczak:2023} for related numerical investigations on the L1 scheme for time discretization. However, to the best of our knowledge, these pointwise-in-time regularity assumptions have not been proved and still await rigorous validation. Moreover, there is no study on the development and analysis of higher-order time-stepping schemes, even for the case with sufficiently smooth problem data.

The present study contributes to the mathematical and numerical analysis of the quasilinear subdiffusion model \eqref{eqn:fde}. First, we  establish new pointwise-in-time regularity estimates on the solution $u$ to problem \eqref{eqn:fde} under suitable regularity assumptions on the problem data (see Assumption \ref{ass:reg-high} and Lemmas \ref{lem: pri-dtu-low reg}-\ref{lem: pri-dttu reg}). This is achieved using a perturbation argument and the $L^p$ regularity estimate in \cite{zacher2012global}. The key of the analysis is the following new H\"{o}lder-type perturbation estimate in Lemma \ref{lem:holder-perturb}: there exists an exponent $\eta\in(0,1)$ such that for any every $\gamma\in[0,1]$, $0\leq s\leq t_*\leq T$, and $q\in(1,\infty)$, there holds
\begin{equation*}
	\|A(u(t_*))^\gamma[ I-A^{-1}(u(t_*))A(u(s)) ] v\|_{L^q(\Omega)} \leq c|t_*-s|^\eta \|A(u(t_*))^\gamma v\|_{L^q(\Omega)},
\end{equation*}
for any $ v\in W^{2,q}(\Omega)\cap H_0^1(\Omega),$ where the notation $A(u(t_*))^\gamma$, $\gamma\geq0$, denotes the fractional power of the elliptic operator $A(u(t_*))$ via spectral decomposition. This estimate allows deriving new pointwise-in-time regularity estimates, including the first- and second time derivatives of the solution. These estimates play a central role in the error analysis of time-stepping schemes. Second, we propose a novel time stepping scheme for the quasilinear model \eqref{eqn:fde} based on the BDF2 convolution quadrature. To restore higher-order convergence, we employ a correction at the starting step \cite{jin2016two,jin2017correction}. We establish a convergence rate $\mathcal{O}(\tau^{1+\alpha-\epsilon})$ for small $\epsilon>0$ in the $L^2(\Omega)$ norm for the time-stepping solution without imposing any additional regularity assumption on the solution $u$, cf. Theorem \ref{thm:main-err-estimate}. The rate is lower than $O(\tau^2)$ due to the limited regularity of the perturbation remainder $R(u;u_0)$, cf. Theorem \ref{thm:reg-rem}. This phenomenon closely resembles that of semilinear subdiffusion \cite[Theorem 3.4]{wangzhou2020high}. Numerical experiments show an empirical rate comparable with the theoretical prediction, confirming the sharpness of the error estimate. The regularity analysis, development of the high-order time-stepping scheme and its convergence analysis represent the main technical contributions of the work.

The rest of the paper is organized  as follows. In Section \ref{sec:sol-reg},
we derive the regularity of the solution $u$ and the remainder $R(u;u_0)$. In Section \ref{sec:err}, we construct a time-stepping scheme via the BDF2 CQ and initial correction and present $L^2(\Omega)$ error bounds. In Section \ref{sec:numers}, we present numerical examples to validate the theoretical results. Throughout, for any $m\geq0$ and $p\geq1$, we denote by $W^{m,p}(\Omega)$ and $W^{m,p}_0(\Omega)$ standard Sobolev spaces of order $m$, equipped with the norm $\|\cdot\|_{W^{m,p}(\Omega)}$ \cite{Adams2003Sobolev}, and  $W^{-m,p'}(\Omega)$ the dual space of $W^{m,p}_0(\Omega)$, with $p'$ being the conjugate exponent of $p$. Further, we write $H^{m}(\Omega)$ and $H^{m}_{0}(\Omega)$ with the norm $\|\cdot\|_{H^m(\Omega)}$ if $p=2$ and write $L^p(\Omega)$ with the norm $\|\cdot\|_{L^p(\Omega)}$ if $m=0$. The notation $(\cdot,\cdot)$ denotes the $L^2(\Omega)$ inner product, and $\|\cdot\|_{\mathcal{B}(L^r)}$ denotes the operator norm on the space $L^r(\Omega)$. We also use Bochner spaces: for a Banach space $B$, let
\begin{equation*}
	W^{m,p}(0,T;B) = \{v: v(t)\in B\ \mbox{for a.e.}\ t\in(0,T)\ \mbox{and}\ \|v\|_{W^{m,p}(0,T;B)}<\infty \}.
\end{equation*}
The space $L^\infty(0,T;B)$ is defined similarly. Throughout, we denote by $c$, with or without a subscript, a generic constant which may differ at each occurrence but is always independent of the time step size $\tau$.

\section{Solution regularity}\label{sec:sol-reg}
Now we discuss the unique solvability and derive Sobolev regularity of the solution $u$ to problem \eqref{eqn:fde}.

\begin{assumption}\label{ass:reg}
	$a\in \mathcal{A}$, $u_0\in W^{2-\frac{2}{p\alpha},p}(\Omega)\cap H_0^1(\Omega)$ and $f\in L^{p}(0,T;L^p(\Omega))$ with $p>d+2\alpha^{-1}$.
\end{assumption}

Under Assumption \ref{ass:reg}, problem \eqref{eqn:fde} admits a unique solution $u \in W^{\alpha,p}(0,T;L^p(\Omega))\cap L^{p}(0,T;W^{2,p}(\Omega)\cap H_0^1(\Omega))$
\cite[Theorem 1.1]{zacher2012global}. By interpolation theory \cite{Amann:2000} and Sobolev embedding theorem \cite[Theorem 4.12, p. 85]{Adams2003Sobolev}, we deduce that for any index $\kappa\in(\frac{1}{\alpha p},\frac{1}{2}-\frac{d}{2p})$, there holds
\begin{equation}\label{eqn:CW1inf}
u\in W^{\kappa\alpha,p}(0,T;W^{2(1-\kappa),p}(\Omega))\hookrightarrow C^{0,\kappa\alpha-\frac{1}{p}}([0,T];W^{1,\infty}(\Omega)).
\end{equation}
This regularity will play a crucial role in deriving of the H\"{o}lder-type perturbation estimates in Lemmas \ref{lem:holder-perturb} and \ref{lem:sharp-holder-perturb}. The main results include pointwise-in-time Sobolev regularity of the solution and its derivatives, and the Sobolev regularity of the remainder $R(u;u_0)$. These results will play a crucial role in the error analysis in Section \ref{sec:err}.

\subsection{Pointwise-in-time solution regularity}
Now we derive several pointwise-in-time regularity estimates on the solution $u$ and its time derivatives. This requires extra regularity conditions on the problem data.
\begin{assumption}\label{ass:reg-high}
$a\in C^2(\mathbb{R})\cap\mathcal{A}$, $u_0\in W^{2,4}(\Omega)\cap W^{2-\frac{2}{p\alpha},p}(\Omega)\cap H_0^1(\Omega)$, with $p>d+2\alpha^{-1}$, and $f\in C^1([0,T];L^4(\Omega))\cap L^{\infty}(0,T;L^\infty(\Omega))$.
\end{assumption}

We shall employ a perturbation argument. For any fix $t_*\in[0,t]$,  we rewrite problem \eqref{eqn:fde} into
\begin{equation}\label{eqn: re-fde}
	\left\{\begin{aligned}
		\Dal  u-\nabla\cdot(a(u(t_*))\nabla u) & = f + R(u;u(t_*)) , \ &\mbox{in}&\ \Omega\times(0,T), \\
		u&=0, \ &\mbox{on}&\ \partial\Omega\times(0,T), \\
		u(0)&=u_0, \ &\mbox{in}&\ \Omega,
	\end{aligned}
	\right.
\end{equation}
with the perturbation remainder  \begin{align}\label{eqn:remain}
R(u;u(t_*)):=\nabla\!\cdot\![(a(u)-a(u(t_*))\nabla u].
\end{align}
Let $A(u)\equiv A(u(t)): H^2(\Omega)\cap H_0^1(\Omega)\mapsto L^2(\Omega)$ be the $L^2(\Omega)$ realization of the elliptic operator
\begin{equation}\label{eqn: con-ellip-oper}
	A(u(t))v = -\nabla\cdot(a(u(t))\nabla v), \quad\forall v\in H^2(\Omega)\cap H_0^1(\Omega).
\end{equation}
For any time $t_*\in[0,t]$, we abbreviate $A(u(t_*))$ to $A_*$. In view of \eqref{eqn:CW1inf}, $A_*$ has full regularity pickup in $W^{2,\infty}(\Omega)\cap H_0^1(\Omega)$. Using Laplace transform, the solution $u(t)$ to problem \eqref{eqn: re-fde} can be represented by \cite[Section 6.2]{Jin:2021book}
\begin{align}
	u(t) &= F(t;A_*)u_0 + \int_0^t E(t-s;A_*)f(s)\ {\rm d}s + \int_0^t E(t-s;A_*) R(u;u(t_*)) \ {\rm d}s \nonumber\\
   &= F(t;A_*)u_0 + \int_0^t E(t-s;A_*)f(s)\ {\rm d}s + \int_0^t E(t-s;A_*) [ A_*-A(u(s)) ] u(s) \ {\rm d}s ,\label{eqn: sol-rep}
\end{align}
where the solution operators $F(t;A_*)$ and $E(t;A_*)$ are defined respectively by
\begin{equation}\label{eqn:oper-EF}
	F(t;A_*)=\frac{1}{2\pi i}\int_{\Gamma_{\theta,\sigma}}e^{zt}z^{\alpha-1}(z^\alpha+A_*)^{-1}\ \mathrm{d}z
	\quad \mbox{and}\quad
	E(t;A_*)=\frac{1}{2\pi i}\int_{\Gamma_{\theta,\sigma}}e^{zt}(z^\alpha+A_*)^{-1}\ \mathrm{d}z,
\end{equation}
with the contour
$\Gamma_{\theta,\sigma} = \{z\in \mathbb{C}:|z|=\sigma,|\arg(z)|\leq\theta\}\cup\{z\in \mathbb{C}:|\arg(z)|=\theta,|z|\geq\sigma\}\subset \mathbb{C}$,  oriented with an increasing imaginary part. Throughout $\theta\in(\frac{\pi}{2},\pi)$ is fixed so that $z^\alpha\in\Sigma_{\alpha\theta}\subset\Sigma_{\theta}$ for all $z\in\Sigma_{\theta}:=\{ z\in\mathbb{C}\setminus\{0\}:|\arg(z)|\leq\theta\}$.

Now we give two preliminary tools that will be used frequently below. First, in view of the estimate \eqref{eqn:CW1inf}, we have the following uniform $L^q(\Omega)$ bound on Riesz transform (see, e.g., \cite[Theorem 4]{AuscherTchamitchian:2001} and \cite[p. 176 and Theorem C]{Shen:2005}).
\begin{lemma}\label{lem:Riesz}
The operator $\nabla A_*^{-\frac12}$ is uniformly bounded on $L^q(\Omega)$, $1<q<\infty$, and thus
\begin{equation*}
		c'\|A_*^{\frac12}v\|_{L^q(\Omega)}\leq \|\nabla v\|_{L^q(\Omega)} \leq c \|A_*^{\frac12}v\|_{L^q(\Omega)}, \quad \forall v\in W_0^{1,q}(\Omega).
\end{equation*}
\end{lemma}

We often use the following Gronwall's inequality (see, e.g., \cite[p. 188]{henry2006geometric} and \cite[Theorem 4.1]{Jin:2021book}).
\begin{lemma}\label{lem:gronwall}
Suppose $b\geq0$, $\beta>0$ and $a(t)$ is a nonnegative function locally integrable on
$[0,T]$,and suppose $u(t)$ is nonnegative and locally integrable on $[0,T]$ with
$u(t) \leq a(t)+b\int_0^t(t-s)^{\beta-1}u(s)\d s$ for all $t \in [0,T]$. Then
\begin{equation*}
  u(t) \leq a(t) + c_T\int_0^t(t-s)^{\beta-1}a(s)\d s, \quad \forall t\in[0,T].
\end{equation*}
\end{lemma}

Next we describe smoothing  properties of the operators $F(t;A_*)$ and $E(t;A_*)$ in $L^q(\Omega)$ spaces. For any $s\in\mathbb{R}$, we define the fractional power $A_*^s$ of $A_*$ via the spectral decomposition.
\begin{lemma}\label{lem: asymp-sol}
Let  Assumption \ref{ass:reg} hold. Then for any $q\in[1,\infty]$, $\ell=0,1$, and $\gamma\in[0,1]$, there exists  $c>0$ such that	
\begin{equation*}
t^{\ell+\gamma\alpha}\|A_*^\gamma F^{(\ell)}(t;A_*)\|_{\mathcal{B}(L^q)} +
t^{\ell+1-(1-\gamma)\alpha}\|A_*^\gamma E^{(\ell)}(t;A_*)\|_{\mathcal{B}(L^q)}\leq c, \quad \forall  t\in(0,T].
\end{equation*}
\end{lemma}
\begin{proof}
These estimates are derived using the $L^q(\Omega)$-resolvent estimate of $A_*$: for all $q\in[1,\infty]$
\begin{equation}\label{eqn:resolvent}
\|  (z + A_*)^{-1} v  \|_{L^q\II} \le c(1+|z|)^{-1}\|
v \|_{L^q\II},\quad \forall v\in L^q\II.
\end{equation}
For $q=\infty$, the estimate \eqref{eqn:resolvent} was shown in \cite[Theorem 2.1]{bakaev2003maximum} for $v\in C_0(\overline \Omega)$ and then in \cite[Appendix A]{LiMa:SINUM2022}
for $v\in L^\infty(\Omega)$. For $q=1$, by the duality argument, and \eqref{eqn:resolvent} with $q=\infty$, there holds for any $v\in L^1(\Omega)$,
	\begin{align*}
		\|(z+A_*)^{-1}v\|_{L^1(\Omega)} &= \sup_{\|w\|_{L^\infty(\Omega)}=1}((z+A_*)^{-1}v,w)= \sup_{\|w\|_{L^\infty(\Omega)}=1}(v,(z+A_*)^{-1}w) \\ &
		\leq \sup_{\|w\|_{L^\infty(\Omega)}=1}\|v\|_{L^1(\Omega)}\|(z+A_*)^{-1}w\|_{L^\infty(\Omega)}\leq c(1+|z|)^{-1}\|v\|_{L^1(\Omega)}.
\end{align*}
The case $q\in(1,\infty)$ follows from Riesz--Thorin interpolation theorem. Using the resolvent estimate \eqref{eqn:resolvent} and repeating the argument in \cite[Theorem 6.4]{Jin:2021book} give the desired estimate.
\end{proof}

The next lemma gives uniform bounds on $u(t)$ and $\Dal  u(t)$.
\begin{theorem}\label{thm: pri-u-h2reg}
Under Assumption \ref{ass:reg-high}, the following time-regularity estimate holds
\begin{equation*}
	\|u(t)\|_{W^{2,4}(\Omega)} + \|\Dal  u(t)\|_{L^4(\Omega)} \leq c, \quad\forall t\in(0,T],
	\end{equation*}
\end{theorem}
\begin{proof}
By the definition of $A(u(s))$, we have the following splitting
\begin{align}
		& [ A_*-A(u(s)) ] u(s) = \nabla [a(u(s))-a(u(t_*))]\cdot\nabla u(s) + [a(u(s))-a(u(t_*))]\Delta u(s)\nonumber \\
		 =& [a'(u(s))\nabla u(s) - a'(u(t_*))\nabla u (t_*)]\cdot\nabla u(s) + [a(u(s))-a(u(t_*))]\Delta u(s)\nonumber \\
		 =& [a'(u(s))\nabla (u(s)-u(t_*))]\cdot\nabla u(s) + [ (a'(u(s))-a'(u(t_*)))\nabla u (t_*)]\cdot\nabla u(s)\nonumber \\
        &+ [a(u(s))-a(u(t_*))]\Delta u(s).\label{eqn: split-nonlinearity}
\end{align}
By Assumption \ref{ass:reg-high}, the regularity estimate $u\in C^{0,\kappa\alpha-\frac{1}{p}}([0,T];W^{1,\infty}(\Omega))$ from \eqref{eqn:CW1inf} and the full elliptic regularity pickup of $A_*$ in $W^{2,\infty}(\Omega)\cap H^1_0(\Omega)$, we obtain
\begin{align}
	&\|[ A_*-A(u(s))] u(s)\|_{L^4(\Omega)}\nonumber\\
  \leq& c(\|\nabla(u(s)-u(t_*))\cdot\nabla u(s)\|_{L^4(\Omega)} + \|(u(s)-u(t_*))\nabla u(s)\|_{L^4(\Omega)} \nonumber\\
  &+ \|(u(s)-u(t_*))\Delta u(s)\|_{L^4(\Omega)})\nonumber\\
		\leq& c\|u(s)-u(t_*)\|_{W^{1,\infty}(\Omega)}\|u(s)\|_{W^{2,4}(\Omega)}
		\leq c|s-t_*|^{\kappa\alpha-\frac{1}{p}}\|A_*u(s)\|_{L^4(\Omega)}.\label{ineq: Holder prop}
\end{align}
Then the identity
$\frac{{\rm d}}{{\rm d}t} F(t;A_*)=-A_*E(t;A_*)$ \cite[Lemma 6.2]{Jin:2021book}, integration by parts, Assumption \ref{ass:reg-high}, and the estimate \eqref{ineq: Holder prop} imply that for any $t_*\in[0,T]$,
\begin{align*}
\|A_*u(t_*)\|_{L^4(\Omega)}&\leq \|F(t_*;A_*)\|_{\mathcal{B}(L^4)}\|A_*u_0\|_{L^4(\Omega)} + \bigg\| \int_0^{t_*}A_*E(t_*-s;A_*)f(s) \ {\rm d}s\bigg\|_{L^4(\Omega)}\\
 &\quad + \int_0^{t_*} \|A_*E(t_*-s;A_*)\|_{\mathcal{B}(L^4)}\|[ A_*-A(u(s))] u(s)\|_{L^4(\Omega)} \ {\rm d}s\\
&\leq c_T +  c\int_0^{t_*}(t_*-s)^{\kappa\alpha-\frac{1}{p}-1}\|A_*u(s)\|_{L^4(\Omega)}\ {\rm d}s.
\end{align*}
By Gr{o}nwall's inequality in Lemma \ref{lem:gronwall} and the full regularity pickup of $A_*$ in $W^{2,\infty}(\Omega)\cap H^1_0(\Omega)$, the first assertion follows.
Moreover, since  $\|\nabla a(u)\|_{L^\infty(\Omega)} = \|a'(u)\nabla u\|_{L^\infty(\Omega)}\leq c$, we deduce
\begin{align*}
\|\Dal  u(t)\|_{L^4(\Omega)} &\leq \|a(u(t))\Delta u(t)\|_{L^4(\Omega)} + \|\nabla a(u)\cdot\nabla u(t)\|_{L^4(\Omega)} + \|f(t)\|_{L^4(\Omega)} \\
&\le c( \| u(t) \|_{W^{2,4}\II} + \|f(t) \|_{L^4(\Omega)} )\leq c,\quad \forall t\in(0,T].
\end{align*}
This completes the proof of the lemma.
\end{proof}

Next we derive a crucial H\"{o}lder type perturbation estimate. Due to the nonlinearity of the diffusion coefficent $a(u)$, the continuity modulus in time is only H\"{o}lderian, which is lower than the Lipschitz continuity for linear problems with time-dependent coefficients \cite[Corollary 9.1]{JinZhou:2023book}.
\begin{lemma}\label{lem:holder-perturb}
Let Assumption \ref{ass:reg-high} hold. For every $\gamma\in[0,1]$ and $q\in(1,\infty)$, there holds for any $0\leq s \leq t_*\leq T$:
	\begin{equation}\label{ineq:holder-perturb}
		\|A_*^\gamma[ I-A_*^{-1}A(u(s)) ] v\|_{L^q(\Omega)} \leq c|t_*-s|^{\kappa\alpha-\frac1p}\|A_*^\gamma v\|_{L^q(\Omega)}, \quad\forall v\in W^{2,q}(\Omega)\cap H_0^1(\Omega).
	\end{equation}
\end{lemma}
\begin{proof}
Let $q$ and $q'$ be a H\"{o}lder conjugate pair. The case $\gamma=1$ follows by repeating the argument for the estimate \eqref{ineq: Holder prop}. For the case $\gamma=0$, let $\psi=A(u(s)) v$ and $\varphi=A_*^{-1}\psi$.
Then we have
\begin{align*}
    \|\varphi-v\|_{L^q(\Omega)} &= \sup_{\|\chi\|_{L^{q'}(\Omega)}=1}(\varphi-v,\chi) = \sup_{\|\chi\|_{L^{q'}(\Omega)}=1} (A_*(\varphi-v),A_*^{-1}\chi) \\
    &= \sup_{\|\chi\|_{L^{q'}(\Omega)}=1}([A(u(s))-A_*]A_*^{-1}\chi,v).
\end{align*}
Repeating the argument for the estimate \eqref{ineq: Holder prop} and the elliptic regularity pickup of $A_*$ in $W^{2,\infty}(\Omega)\cap H_0^1(\Omega)$ imply
\begin{equation*}
	\|[ A_*-A(u(s))]A_*^{-1}\chi\|_{L^{q'}(\Omega)}\leq c|t_*-s|^{\kappa\alpha-\frac1p}\|A_*A_*^{-1}\chi\|_{L^{q'}(\Omega)}= c|t_*-s|^{\kappa\alpha-\frac1p}\|\chi\|_{L^{q'}(\Omega)}.
\end{equation*}
Then H\"{o}lder's inequality leads to
\begin{equation*}
    	\|\varphi-v\|_{L^q(\Omega)}\leq c|t_*-s|^{\kappa\alpha-\frac1p}\|v\|_{L^q(\Omega)}.
 \end{equation*}
This proves the  estimate for the case $\gamma=0$. The case  $\gamma\in(0,1)$ follows from interpolation.
\end{proof}

Next we bound the time derivatives $u'(t)$ and $u''(t)$. We use frequently the following shorthand
\begin{equation}\label{eqn:hatA}
  \widehat A(s) = -\frac{\d}{\d s}A(u(s))\quad \mbox{and}\quad \widetilde A(s) = -\frac{\d^2}{\d s^2}A(u(s)).
\end{equation}
Then direct computation gives
\begin{align}
\widehat{A}(s) u(s) = [a''(u(s))|\nabla u(s)|^2+a'(u(s))\Delta u(s)] u'(s) + a'(u(s))\nabla u(s)\cdot\nabla  u'(s). \label{eqn:oper-expansion}
\end{align}
In the next lemma, we give several preliminary estimates on $\widehat{A}(t)$ and $\widetilde{A}(t)$.
\begin{lemma} \label{lem: non-perb estimate}
Let $\Omega$ be a $C^2$ domain and Assumption \ref{ass:reg-high} hold, $\gamma\in[0,1/2)$, and $1<r<\infty$. Then there hold for $\ell=0,1$ and $0<s<t_*\leq T$:
    \begin{enumerate}
        \item[{\rm(i)}] $\|A_*^\gamma E^{(\ell)}(t_*-s;A_*)\widehat A(s) u(s)\|_{L^r(\Omega)} \leq c(t_*-s)^{(\frac12-\gamma)\alpha-1-\ell}\|u'(s)\|_{L^r(\Omega)}$;
        \item[{\rm(ii)}] $\|A_*^\gamma E(t_*-s;A_*)\widetilde {A}(s)u(s)\|_{L^r(\Omega)} \leq c(t_*-s)^{(\frac{1}{2}-\gamma)\alpha-1}(\|u'(s)\|_{L^{2r}(\Omega)}^2+\|u''(s)\|_{L^r(\Omega)})$;
        \item[{\rm(iii)}] $\|A_*^\gamma E(t_*-s;A_*)\widehat{A}(s)u'(s)\|_{L^r(\Omega)} \leq c(t_*-s)^{(\frac12-\gamma)\alpha-1}
        \|u'(s)\|_{L^{2r}(\Omega)}\|\nabla u'(s)\|_{L^{2r}(\Omega)}$.
\end{enumerate}
\end{lemma}
\begin{proof}
The estimates follows from the smoothing properties of $E(t;A_*)$ and a duality argument. Indeed, by Lemma \ref{lem: asymp-sol}, H\"{o}lder's inequality, Lemma \ref{lem:Riesz}, the regularity estimate \eqref{eqn:CW1inf}
and Assumption \ref{ass:reg-high}, we get
\begin{align}
  &\big\|A_*^\gamma E^{(\ell)}(t_*-s;A_*)\widehat{A}(s) u(s)\big\|_{L^r(\Omega)}\nonumber\\
  =& \sup_{\|\varphi\|_{L^{r'}(\Omega)=1}} (a'(u(s))u'(s)\nabla u(s)  ,\nabla A_*^\gamma E^{(\ell)}(t_*-s;A_*)\varphi) \nonumber\\
     \leq&c\|\nabla u(s)\|_{L^\infty(\Omega)}\|u'(s)\|_{L^r(\Omega)}\|A_*^{\frac12+\gamma}E^{(\ell)}(t_*-s;A_*)\|_{\mathcal{B}(L^{r'})}\nonumber\\
     \leq& c(t_*-s)^{(\frac12-\gamma)\alpha-1}\|u'(s)\|_{L^r(\Omega)}.\label{ineq: useful estimate}
\end{align}
Similarly, by Lemma \ref{lem: asymp-sol}, Assumption \ref{ass:reg-high}, and H\"{o}lder's inequality, we deduce
\begin{align*}
	&\|A_*^\gamma E(t_*-s;A_*)\widetilde{A}(s)u(s)\|_{L^r(\Omega)}\nonumber\\
     = &\sup_{\|\varphi\|_{L^{r'}(\Omega)=1}} ([a''(u(s))|u'(s)|^2+a'(u(s))u''(s)]\nabla u(s) ,\nabla A_*^\gamma E(t_*-s;A_*)\varphi) \nonumber\\
	 \leq &c(\|u'(s)\|^2_{L^{2r}(\Omega)} + \|u''(s)\|_{L^r(\Omega)})\|\nabla u(s)\|_{L^\infty(\Omega)}
\|A_*^{\frac12+\gamma} E(t_*-s;A_*)\|_{\mathcal{B}(L^{r'})}. 
\end{align*}
The bound on $\|A_*^\gamma E(t_*-s;A_*)\widehat{A}(s) u'(s)\|_{L^r(\Omega)}$ follows analogously, and hence the proof is omitted.
\end{proof}

Now we bound the first time derivative $u'(t)$. The proof employs Lemma \ref{lem: asymp-sol} and a perturbation argument.
\begin{lemma}\label{lem: pri-dtu-low reg}
Under Assumption \ref{ass:reg-high}, there exists $c=c(T,\alpha,\Omega,p,\kappa,f,u_0,\gamma)$ such that for $0<t\leq T$ and $0\leq \gamma<1$:
	\begin{equation*}
	    \|A_*^\gamma u'(t)\|_{L^4(\Omega)} \leq c\max(t^{(1-\gamma)\alpha-1},1).
	\end{equation*}
\end{lemma}
\begin{proof}
Fix  $r=4$ and $r'=4/3$. It suffices to prove the cases
$\gamma=0$ and $\gamma_0\in(1/2,1)$, and the case $\gamma\in(0,\gamma_0)$ follows by interpolation. Differentiating the identity \eqref{eqn: sol-rep} with respect to $t$ gives
\begin{equation}\label{eqn: dtu-repre}
    \begin{split}
      u'(t) =& -A_*E(t;A_*)u_0 + E(t;A_*)f(0) + \int_0^t E(t-s;A_*)f'(s)\ {\rm d}s \nonumber\\&\quad+ E(t;A_*)( A_* - A(u_0)) u_0 + \int_0^t E(t-s;A_*) \frac{{\rm d}}{{\rm d}s}\{[ A_*-A(u(s)) ] u(s)\} \ {\rm d}s
	\\=& \Big[E(t;A_*)(f(0)-A(u_0)u_0) + \int_0^t E(t-s;A_*)f'(s)\ {\rm d}s\Big] \nonumber\\
    &+ \int_0^t E(t-s;A_*) [ A_*-A(u(s)) ] u'(s)\ {\rm d}s  \nonumber\\&+ \int_0^t E(t-s;A_*) \widehat{A}(s) u(s) \ {\rm d}s := {\rm I}(t) + {\rm II}(t) + {\rm III}(t).
    \end{split}
\end{equation}
Below we fix $\gamma\in\{0,\gamma_0\}$. It follows from Lemma \ref{lem: asymp-sol} and Assumption \ref{ass:reg-high} that
\begin{align*}
	\|A_*^\gamma {\rm I} (t_*)\|_{L^r(\Omega)}
\leq& c\|A_*^\gamma E(t_*;A_*)\|_{\mathcal{B}(L^r)} \|f(0)-A(u_0)u_0\|_{L^r(\Omega)}\\
  &+ c\int_0^{t_*} \|A^\gamma_*E(t_*-s;A_*)\|_{\mathcal{B}(L^r)} \|f'(s)\|_{L^r(\Omega)}\ {\rm d}s \\
 \leq& c(t_*^{(1-\gamma)\alpha-1}+c_T)\leq c_T(t_*^{(1-\gamma)\alpha-1}+1).
\end{align*}
By Lemma \ref{lem:holder-perturb}, we can bound the term ${\rm II}(t_*)$ by
\begin{align*}
		\|A_*^\gamma{\rm II}(t_*)\|_{L^r(\Omega)}&\leq \int_0^{t_*} \|A_*E(t_*-s;A_*)\|_{\mathcal{B}(L^r)} \|A_*^\gamma[ I-A_*^{-1}A(u(s)) ] u'(s)\|_{L^r(\Omega)}\ {\rm d}s  \\
		&\leq c \int_0^{t_*}(t_*-s)^{\kappa\alpha-\frac1p-1}\|A_*^\gamma u'(s)\|_{L^r(\Omega)} \ {\rm d}s.
	\end{align*}
Next, we treat the cases $\gamma=0$ and $\gamma=\gamma_0\in(1/2,1)$ separately.
When $\gamma=0$, Lemma \ref{lem: non-perb estimate}(i) with $r$ directly implies
\begin{equation*}
   \|{\rm III}(t_*)\|_{L^r(\Omega)}\leq c\int_0^{t_*}(t_*-s)^{\frac{\alpha}{2}-1}\|u'(s)\|_{L^r(\Omega)} \ {\rm d}s.
\end{equation*}
When $\gamma=\gamma_0>\frac12$, using the identity \eqref{eqn:oper-expansion},
by Assumption \ref{ass:reg-high}, we have the regularity $u\in C^{0,\kappa\alpha-\frac{1}{p}}([0,T];W^{1,\infty}(\Omega))$ from \eqref{eqn:CW1inf}. This, Theorem \ref{thm: pri-u-h2reg}, Lemma \ref{lem:Riesz} and Sobolev embedding $W^{1,r}(\Omega)\hookrightarrow L^{\infty}(\Omega)$ \cite[Theorem 4.12, p. 85]{Adams2003Sobolev}, we obtain
\begin{align*}
\|\widehat{A}(s) u(s)\|_{L^r(\Omega)} & \leq c(1+\|\Delta u(s)\|_{L^r(\Omega)})\|u'(s)\|_{L^{\infty}(\Omega)} + c\|\nabla u'(s)\|_{L^r(\Omega)}\\
&\leq c\|A_*^\frac12u'(s)\|_{L^r(\Omega)}\leq c\|A_*^{\gamma_0} u'(s)\|_{L^r(\Omega)},
\end{align*}
by the continuous embedding $(D(A_*^{\gamma_0}),L^r(\Omega))\hookrightarrow  (D(A_*^\frac12).L^r(\Omega))$ \cite{Sobolevskii:1964}.
Consequently,
\begin{align*}
	\|A_*^\gamma {\rm III}(t_*)\|_{L^r(\Omega)} &\leq \int_0^{t_*}\|A_*^\gamma E(t_*-s;A_*)\|_{\mathcal{B}(L^r)}\|\widehat{A}(s) u(s)\|_{L^r(\Omega)}\ {\rm d}s \\
 & \leq c\int_0^{t_*}(t_*-s)^{(1-\gamma)\alpha-1}\|A_*^\gamma u'(s)\|_{L^r(\Omega)}\ {\rm d}s.
 \end{align*}
Then for $\gamma\in\{0,\gamma_0\}$, Gr\"{o}nwall's inequality yields the desired assertion.
\end{proof}

Next we establish a sharp H\"{o}lder type perturbation estimate, with the continuity modulus in time depending only on the order $\alpha$, removing the dependence of the constant $c$ on $\kappa$ and $p$.
\begin{lemma}\label{lem:sharp-holder-perturb}
Let Assumption \ref{ass:reg-high} hold. For every $\gamma\in[0,1]$ and $q\in[2,4]$, there holds for any $0\leq s \leq t_*\leq T$:
	\begin{equation}\label{ineq:sharp-holder-perturb}
		\|A_*^\gamma[ I-A_*^{-1}A(u(s)) ] v\|_{L^q(\Omega)} \leq c|t_*-s|^{\frac{\alpha}{2}}\|A_*^\gamma v\|_{L^q(\Omega)}, \quad\forall v\in W^{2,q}(\Omega)\cap H_0^1(\Omega).
	\end{equation}
\end{lemma}
\begin{proof}
By interpolation, it suffices to prove the assertion in the $L^2(\Omega)$ and $L^4(\Omega)$ norms. By Lemma \ref{lem: pri-dtu-low reg}, Sobolev embedding  $W^{1,4}\II\hookrightarrow L^\infty\II$ (for $d=1,2,3$) and Poinc{a}r\'{e} inequality, we have			
\begin{align}
     \|\nabla(u(s)-u(t_*))\|_{L^4\II} &\leq \int_s^{t_*} \|\nabla u'(\xi)\|_{L^4\II}\ {\rm d}\xi \leq c\Big|\int_s^{t_*} \xi^{\frac{\alpha}{2}-1}\ {\rm d}\xi\Big| \leq c|s-t_*|^{\frac{\alpha}{2}},\label{eqn:basic-est-1} \\
    \|u(s)-u(t_*)\|_{L^\infty\II} &\leq c\|\nabla(u(s)-u(t_*))\|_{L^4\II} \leq c|s-t_*|^{\frac{\alpha}{2}}. \label{eqn:basic-est-2}
\end{align}
Meanwhile, we have (cf. \eqref{eqn: split-nonlinearity})
\begin{align*}
	&[ A_*-A(u(s)) ] v
	 = [a'(u(s))\nabla (u(s)-u(t_*))]\cdot\nabla v \\
    &+ ([ (a'(u(s))-a'(u(t_*)))\nabla u (t_*)]\cdot\nabla v + [a(u(s))-a(u(t_*))]\Delta v) :={\rm I} + {\rm II}.
\end{align*}
Now Assumption \ref{ass:reg-high}, H\"{o}lder's inequality, the estimates \eqref{eqn:basic-est-1}--\eqref{eqn:basic-est-2}, and Sobolev embedding $ W^{1,q}(\Omega)\hookrightarrow L^r(\Omega)$ (with $r=4$ for $q=2$, and $r=\infty$ for $q=4$) \cite[Theorem 4.12, p. 85]{Adams2003Sobolev} imply that for $q\in \{2,4\}$,
\begin{align*}
    \|{\rm I}\|_{L^q\II}
	&\leq c\|\nabla(u(s)-u(t_*))\|_{L^4\II}\| \nabla v\|_{L^r\II} \leq c|s-t_*|^{\frac{\alpha}{2}}\|A_*v\|_{L^q\II},\\
\|{\rm II}\|_{L^q\II}
	&\leq c\|u(s)-u(t_*)\|_{L^\infty\II}(\| \nabla v\|_{L^q\II} + \| \Delta v\|_{L^q\II}) \leq c|s-t_*|^{\frac{\alpha}{2}}\|A_*v\|_{L^q\II}.
\end{align*}
This proves the assertion for $\gamma=1$. For $\gamma=0$, by a duality argument, with $\psi=A(u(s)) v$ and $\varphi=A_*^{-1}\psi$,
\begin{equation*}
    \|\varphi-v\|_{L^q\II}= \sup_{\|\chi\|_{L^{q'}\II}=1}([A(u(s))-A_*]A_*^{-1}\chi,v).
\end{equation*}
Then it suffices to bound $\rm I$ and $\rm II$ with $v=A_*^{-1}\chi$ in $L^{q'}(\Omega)$.
Assumption \ref{ass:reg-high}, H\"{o}lder's inequality, the estimates \eqref{eqn:basic-est-1}--\eqref{eqn:basic-est-2}, and Sobolev embedding
$W^{1,q'}(\Omega)\hookrightarrow L^{\tilde r}(\Omega)$ (with $\tilde r=2$ for $q=4$, and $\tilde r=4$ for $q=2$) lead to
\begin{align*}
    \|{\rm I}\|_{L^{q'}\II}
	&\leq c\|\nabla(u(s)-u(t_*))\|_{L^4\II}\| \nabla A_*^{-1}\chi\|_{L^{\tilde r}\II} \\
 &\leq c|s-t_*|^{\frac{\alpha}{2}}\|A_*^{-1}\chi\|_{W^{2,q'}\II} \leq c|s-t_*|^{\frac{\alpha}{2}}\|\chi\|_{L^{q'}\II},\\
    \|{\rm II}\|_{L^{q'}\II}
	 &\leq c\|u(s)-u(t_*)\|_{L^\infty\II}(\| \nabla A_*^{-1}\chi\|_{L^{q'}\II} + \| \Delta A_*^{-1}\chi\|_{L^{q'}\II}) \\& \leq c|s-t_*|^{\frac{\alpha}{2}}\|A_*^{-1}\chi\|_{W^{2,q'}\II} \leq c|s-t_*|^{\frac{\alpha}{2}}\|\chi\|_{L^{q'}\II}.
\end{align*}
Thus we have proved the case $\gamma=0$. The cases with $\gamma\in(0,1)$ and $q\in(2,4)$ follow by interpolation.
\end{proof}

Last, we derive pointwise-in-time regularity of $u''(t)$ under the following assumption.
\begin{assumption}\label{assum: reg-plus}
Assumption \ref{ass:reg-high} holds, and $f\in C^2([0,T];L^2(\Omega))$.
\end{assumption}

Now we can bound $u''$. The main challenge lies in the strong singularity of $E'(t;A)$ at $t=0$.
\begin{lemma}\label{lem: pri-dttu reg}
Fix $\epsilon_0>0$ small. Let $\Omega$ be a $C^2$ domain, and Assumption \ref{assum: reg-plus} hold. There exists $c=c(T,\alpha,\Omega,f,u_0,\gamma)>0$ such that for $0<t\leq T$ and $0\leq\gamma\leq\frac12-\epsilon_0$:
	\begin{equation*}
		\|A_*^\gamma u''(t)\|_{L^2(\Omega)} \leq c\max(t^{(1-\gamma)\alpha-2},1).
	\end{equation*}
\end{lemma}
\begin{proof}
By interpolation, it suffices to prove the estimate for $\gamma\in\{0,\frac12-\epsilon_0\}$. Let $B_\gamma(t) = A^\gamma_* E(t;A_*)$,  $\widetilde \gamma = 1-(1-\gamma)\alpha\in(01)$ and $\widehat\gamma=1-(\frac{1}{2}-\gamma)\alpha$, and we use $\widehat{A}(t)$ and $\widetilde{A}(t)$ defined in \eqref{eqn:hatA}.
By Lemma \ref{lem: asymp-sol},
\begin{equation}\label{eqn:B-bdd}
\|B_\gamma(t)\|_{\mathcal{B}(L^2)}\leq ct^{-\widetilde\gamma} \quad \mbox{and}\quad \|B'_\gamma(t)\|_{\mathcal{B}(L^2)}\leq ct^{-\widetilde\gamma-1}.
\end{equation}
Meanwhile, from Lemmas \ref{lem: non-perb estimate} and \ref{lem: pri-dtu-low reg}, we have the following basic estimates
\begin{align}
\|B_\gamma^{(\ell)}(t_*-s)\widehat A(s) u(s)\|_{L^2(\Omega)}&\leq c(t_*-s)^{-\widehat\gamma -\ell}\max(s^{\alpha-1},1),\label{eqn:basi-est1}\\
\|B_\gamma(t_*-s)\widetilde {A}(s)u(s)\|_{L^2(\Omega)}&\leq c(t_*-s)^{-\widehat\gamma}(s^{2\alpha-2}+\|u''\|_{L^2(\Omega)}),\label{eqn:basic-est2}\\
\|B_\gamma(t_*-s)\widehat{A}(s)u'(s)\|_{L^2(\Omega)} &\leq c(t_*-s)^{-\widehat\gamma}
        \max(s^{\frac{3\alpha}{2}-2},1).\label{eqn:basic-est3}
\end{align}
Differentiating \eqref{eqn: sol-rep} in $t$ and multiplying  $t^{\widetilde\gamma+1}$ give
\begin{align}
	t^{\widetilde\gamma+1}u'(t) =& t^{\widetilde\gamma+1}E(t;A_*)(f(0)-A_0u_0) + t^{\widetilde\gamma+1}\int_0^t E(t-s;A_*)f'(s)\ {\rm d}s \nonumber\\
   &+ t^{\widetilde\gamma+1}\int_0^t E(t-s;A_*) \widehat{A}(s) u(s) \ {\rm d}s \nonumber\\
   &+ t^{\widetilde\gamma+1}\int_0^t E(t-s;A_*) [ A_*-A(u(s)) ] u'(s)\ {\rm d}s := {\rm I}(t)+{\rm II}(t)+{\rm III}(t)+{\rm IV}(t).\nonumber 
\end{align}
The estimates in \eqref{eqn:B-bdd} and Assumption \ref{ass:reg-high} imply
\begin{align*}
	\|A_*^\gamma{\rm I}'(t_*)\|_{L^2(\Omega)}\leq& c(t_*^{\widetilde\gamma}\|B_\gamma (t_*)\|_{\mathcal{B}(L^2)}+t_*^{\widetilde\gamma+1} \|B_\gamma '(t_*)\|_{\mathcal{B}(L^2)})\|f(0)-A_0u_0\|_{L^2(\Omega)}\leq c,\\
    \|A_*^\gamma {\rm II}'(t_*)\|_{L^2(\Omega)} \leq &ct_*^{\widetilde\gamma}\int_0^{t_*}\!\!\! \!\|B_\gamma(t_*-s)\|_{\mathcal{B}(L^2)}\|f'(s)\|_{L^2(\Omega)}\ {\rm d}s
   \!+\! ct_*^{\widetilde\gamma+1}\|B_\gamma (t_*)\|_{\mathcal{B}(L^2)} \| f'(0)\|_{L^2(\Omega)} \\
    &+ ct_*^{\widetilde\gamma+1}\int_0^{t_*}\|B_\gamma(t_*-s)\|_{\mathcal{B}(L^2)}\|f''(s)\|_{L^2(\Omega)}\ {\rm d}s \leq c_T.
\end{align*}
For the term ${\rm III}(t)$, with the identity $t = (t-s)+s$, we have
\begin{align*}
	{\rm III}(t)  =& t^{\bar\gamma}\int_0^t (t-s)E(t-s;A_*)\widehat{A}(s)u(s) \ {\rm d}s  + t^{\bar\gamma }\int_0^t E(s;A_*) (t-s)\widehat{A}(t-s) u(t-s) \ {\rm d}s\\
 =:& {\rm III}_{1}(t) + {\rm III}_{2}(t).
\end{align*}
Then the triangle inequality and the estimate \eqref{eqn:basi-est1} yield 	\begin{align*}
	\|A_*^\gamma{\rm III}'_{1}(t_*)\|_{L^2(\Omega)}& \leq \widetilde\gamma t_*^{\widetilde\gamma-1}\int_0^{t_*} (t_*-s)\|B_\gamma(t_*-s) \widehat{A}(s) u(s)\|_{L^2(\Omega)} \ {\rm d}s \\&\quad + t_*^{\widetilde\gamma }\int_0^{t_*}\|(B_\gamma(t_*-s) +(t_*-s)B_\gamma'(t_*-s)) \widehat{A}(s) u(s)\|_{L^2(\Omega)}\ {\rm d}s
		\\ & \leq ct_*^{\widetilde\gamma-1}\int_0^{t_*} (t_*-s)^{-\widehat\gamma+1}s^{\alpha-1} \ {\rm d}s + ct_*^{\widetilde\gamma }\int_0^{t_*} (t_*-s)^{-\widehat\gamma}s^{\alpha-1} \ {\rm d}s \leq c_T.
\end{align*}
Meanwhile, the estimates \eqref{eqn:basic-est2} and \eqref{eqn:basic-est3} imply
\begin{align*}
	&\|A_*^\gamma {\rm III}'_{2}(t_*)\|_{L^2(\Omega)} \leq\widetilde\gamma t_*^{\widetilde\gamma-1}\int_0^{t_*}s\big\|B_\gamma(t_*-s)\widehat{A}(s) u(s)\big\|_{L^2(\Omega)} \ {\rm d}s\\
  & \quad + t_*^{\widetilde\gamma}\int_0^{t_*}\|B_\gamma(t_*-s) \widehat{A}(s) u(s)+sB_\gamma(t_*-s) (\widehat{A}(s) u'(s)+\widetilde{A}(s) u(s))\|_{L^2(\Omega)} \ {\rm d}s\\
  \leq&  ct_*^{\widetilde\gamma-1}\int_0^{t_*}(t_*-s)^{-\widehat\gamma}s^{\alpha} \ {\rm d}s
  + ct_*^{\widetilde\gamma }\int_0^{t_*}(t_*-s)^{-\widehat\gamma}(s^{\alpha-1}+s\|u''(s)\|_{L^2(\Omega)})\ {\rm d}s\\
	\leq &c_{T} + ct_*^{\widetilde\gamma }\int_0^{t_*}(t_*-s)^{-\widehat\gamma} s\|u''(s)\|_{L^2(\Omega)}\ {\rm d}s.
\end{align*}
For the last term ${\rm IV}$, we split it into
\begin{align*}
		{\rm IV}(t) & = t^{1-\alpha}\int_0^t (t-s)E(t-s;A_*) [ A_*-A(u(s)) ] u'(s) \ {\rm d}s \\&\quad + t^{1-\alpha}\int_0^t E(s;A_*) (t-s)[ A_*-A(u(t-s)) ] u'(t-s) \ {\rm d}s =: {\rm IV}_{1}(t) + {\rm IV}_{2}(t).
	\end{align*}
Repeating the proceeding argument gives
\begin{align*}
   &\|A_*^\gamma{\rm IV}'_{1}(t_*)\|_{L^2(\Omega)}\leq \widetilde\gamma t_*^{\widetilde\gamma-1}\int_0^{t_*} (t_*-s)\|B_\gamma(t_*-s) [ A_*-A(u(s)) ] u'(s)\|_{L^2(\Omega)} \ {\rm d}s\\
   &\quad + t_*^{\widetilde\gamma}\int_0^{t_*} \|(B_\gamma (t_*-s) +(t_*-s)B_\gamma '(t_*-s)) [ A_*-A(u(s)) ] u'(s)\|_{L^2(\Omega)} \ {\rm d}s\\
  &\leq ct_*^{\widetilde\gamma-1}\int_0^{t_*} (t_*-s)^{-\widehat\gamma+1}s^{\alpha-1} \ {\rm d}s + ct_*^{\widetilde\gamma }\int_0^{t_*} (t_*-s)^{-\widehat\gamma}s^{\alpha-1} \ {\rm d}s \leq c_T.
\end{align*}
Last, using the estimates \eqref{eqn:basic-est2} and \eqref{eqn:basic-est3}, we get
\begin{align*}
&\|A_*^\gamma{\rm IV}'_{2}(t_*)\|_{L^2(\Omega)} \leq \widetilde \gamma t_*^{\widetilde\gamma-1}\int_0^{t_*}s\|B_\gamma (t_*-s) [ A_*-A(u(s)) ] u'(s)\|_{L^2(\Omega)} \ {\rm d}s\\
&+ t_*^{\widetilde\gamma}\int_0^{t_*}\!\!\!\!\|B_\gamma(t_*\!-\!s) [ A_*\!-\!A(u(s)) ] u'(s)\!+\!sB_\gamma(t_*\!-\!s) (\widehat{A}(s) u'(s)\!+\! [ A_*\!-\!A(u(s)) ] u''(s))\|_{L^2(\Omega)} {\rm d}s\\
&\leq ct_*^{\widetilde\gamma-1}\int_0^{t_*}(t_*-s)^{-\widehat\gamma}s^{\alpha} \ {\rm d}s
+ ct_*^{\widetilde\gamma}\int_0^{t_*}(t_*-s)^{-\widehat\gamma}(s^{\alpha-1} + s\|u''(s)\|_{L^2(\Omega)}) \ {\rm d}s \\
&\leq c_{T} + ct_*^{\widetilde\gamma }\int_0^{t_*} (t_*-s)^{-\widehat\gamma}[s\|u''(s)\|_{L^2(\Omega)}] \ {\rm d}s.
\end{align*}
Combining the preceding estimates yields
\begin{equation*}
	t_*\|u''(t_*)\|_{L^2(\Omega)} \leq c_{T}t_*^{-\widetilde\gamma} + c\int_0^{t_*} (t_*-s)^{-\widehat\gamma}[s\|u''(s)\|_{L^2(\Omega)}] \ {\rm d}s.
\end{equation*}
The desired assertion now follows directly from Gr\"{o}nwall's inequality in Lemma \ref{lem:gronwall}.
\end{proof}

\subsection{Regularity of the remainder $R(u;u_0)$}

Now we derive the regularity of the remainder $R(u;u_0)$ in Bochner-Sobolev spaces. Given a Banach space $B$, for any $s\ge 0$ and
$1\le p < \infty$, we denote by $W^{s,p}(0,T;B)$ the space of functions $v:(0,T)\rightarrow B$,
with the norm defined by interpolation. For any $0<s< 1$, the Sobolev--Slobodeckij seminorm $|\cdot|_{W^{s,p}(0,T;B)}$ is defined by
\begin{equation}\label{eqn:SS-seminorm}
	| v  |_{W^{s,p}(0,T;B)}^p := \int_0^T\hskip-5pt\int_0^T \frac{\|v(t)-v(\xi)\|_{B}^p}{|t-\xi|^{1+ps}}\ {\rm d} t{\rm d}\xi ,
\end{equation}
and the full norm $\|\cdot\|_{W^{k+s,p}(0,T;B)}$, with $k\ge 0$ and $k\in \mathbb{N}$, is defined by
\begin{equation*}
	\|v\|_{W^{k+s,p}(0,T;B)}^p = \sum_{m=0}^k\|\partial_t^m v \|_{L^p(0,T;B)}^p+|\partial_t^k v |_{W^{s,p}(0,T;B)}^p .
\end{equation*}

Under Assumption \ref{assum: reg-plus}, we have the following regularity estimates on the remainder $R(u;u_0)$ in negative Sobolev spaces. These estimates play a crucial role in the analysis of the corrected BDF2 scheme and also indicate that one can expect a convergence rate $\mathcal{O}(\tau^{1+\alpha-\epsilon})$ for small $\epsilon>0$ at best, cf. Theorem \ref{thm:main-err-estimate} for the precise error estimate. Here for any $s\in\mathbb{R}$, we define the norm $\|v\|_{\dot H^s(\Omega)}=\|A_*^\frac{s}{2}v\|_{L^2(\Omega)}$ (i.e., via spectral decomposition). The restriction $\mu>0$ in Theorem \ref{thm:reg-rem} is due to the limited regularity of the second derivative $u''$ in Lemma \ref{lem: pri-dttu reg}.

\begin{theorem}\label{thm:reg-rem}
Let $\epsilon$ and $\mu>0$ be small. Let Assumption \ref{assum: reg-plus} hold and $u$ be the solution to problem \eqref{eqn:fde}. Then for the remainder $R(u;u_0)$ defined in \eqref{eqn:remain}, there holds
\begin{equation*}
\min(t^{1-\alpha},1)\Big\|\frac{{\rm d}}{{\rm d}t}R(u;u_0)\Big\|_{\dot{H}^{-1}(\Omega)} + \min(t^{2-\alpha},1)\Big\|\frac{{\rm d}^2}{{\rm d}t^2}R(u;u_0)\Big\|_{\dot{H}^{-1-\mu}(\Omega)} \leq c, \quad\forall t\in(0,T].
\end{equation*}
Moreover,  $R(u;u_0)\in W^{1+\alpha-\epsilon,1}(0,T;\dot{H}^{-1-\mu}(\Omega))$.
\end{theorem}
\begin{proof}
Let
$w(t):=\frac{{\rm d}}{{\rm d}t}R (u(t);u_0)= \widehat{A}(t)u(t) + [A(u_0)-A(u(t))]u'(t).$ H\"{o}lder's inequality and Lemma \ref{lem:Riesz} yield
\begin{align*}
\|A(u_0)^{-\frac12}\widehat{A}(t)u(t)\|_{L^2(\Omega)} &= \sup_{\|\varphi\|_{L^2(\Omega)=1}} (a'(u(t))u'(t)\nabla u(t) ,\nabla A(u_0)^{-\frac12}\varphi)\\ &\leq c\|u'(t)\|_{L^2(\Omega)}\leq c\max(t^{\alpha-1},1).
\end{align*}
From Lemmas  \ref{lem: pri-dtu-low reg} and \ref{lem:sharp-holder-perturb}, we deduce
\begin{equation*}
\|A(u_0)^{-\frac12}[A(u_0)-A(u)]u'(t)\|_{L^2(\Omega)}\leq c\max(t^{\alpha-1},1).
\end{equation*}
This proves the desired estimate on $\|w(t)\|_{\dot{H}^{-1}(\Omega)}$. By the product rule, we have
\begin{equation*}
	w'(t)= \widetilde{A}(t)u(t) +2\widehat{A}(t)u' + [A(u_0) -A(u(t))]u''(t).
\end{equation*}
By H\"{o}lder's inequality, Lemma \ref{lem:Riesz} (with $p=2$), Lemmas \ref{lem: pri-dtu-low reg} and \ref{lem: pri-dttu reg} and the regularity estimate \eqref{eqn:CW1inf} yield
\begin{align*}
    	\|A(u_0)^{-\frac12}\widetilde{A}(t)]u(t)\|_{L^2(\Omega)} &= \sup_{\|\varphi\|_{L^2(\Omega)=1}} ([a''(u(t))|u'(t)|^2+a'(u(t))u''(t)]\nabla u(t) ,\nabla A(u_0)^{-\frac12}\varphi) \\
    	& \leq c(\|u'(t)\|^2_{L^4(\Omega)} + \|u''(t)\|_{L^2(\Omega)})
    	\leq c\max(t^{\alpha-2},1),
    \end{align*}
Repeating the argument with Lemma \ref{lem: pri-dtu-low reg} give
    \begin{align*}
        \|A(u_0)^{-\frac12}\widehat{A}(t)]u'(t)\|_{L^2(\Omega)} & = \sup_{\|\varphi\|_{L^2(\Omega)=1}} (a'(u(t))u'(t)\nabla u'(t) ,\nabla A(u_0)^{-\frac12}\varphi)  \\& \leq c\|u'(t)\|_{L^4(\Omega)}\|\nabla u'(t)\|_{L^4(\Omega)}
    	\leq c\max(t^{\frac32\alpha-2},1).
    \end{align*}
Next applying Lemmas \ref{lem:sharp-holder-perturb} and \ref{lem: pri-dttu reg} leads to
\begin{align*}
    	\|A(u_0)^{-\frac12-\frac{\mu}{2}}[A(u_0) -A(u(t))]u''(t)\|_{L^2(\Omega)} &\leq ct^{\frac{\alpha}{2}}\|A(u_0)^{\frac12-\frac{\mu}{2}}u''(t)\|_{L^2(\Omega)}\\
     &\leq c\max(t^{(1+\frac{\mu}{2})\alpha-2},1).
    \end{align*}
Combining the preceding estimates yields
\begin{equation*}
   \|A(u_0)^{-\frac12-\frac{\mu}{2}}w'(t)\|_{L^2(\Omega)} \leq c\max(t^{\alpha-2},1).
\end{equation*}
This proves the first assertion. Next we bound $R(u;u_0)$ in the $W^{1+\alpha-\epsilon,1}(0,T;\dot{H}^{-1-\mu}(\Omega))$ norm. By integration by parts, we have the following estimate
\begin{align*}
	&\int_{0}^{1} \frac{t^{\alpha-1} -1 }{(1-t)^{1+\alpha-\epsilon}} \ {\rm d} t
    \le  \Big( \int_{0}^{\frac12}+  \int_{\frac12}^{1} \Big) \frac{ t^{\alpha-1} -1 }{(1-t)^{1+\alpha-\epsilon}} \ {\rm d} t\\
    \le& c + c \lim_{t\rightarrow 1^-} \frac{t^{\alpha-1} -1 }{(1-t)^{\alpha-\epsilon}} + c\int_{\frac12}^{1} \frac{ t^{\alpha-2} }{(1-t)^{\alpha-\epsilon}} \ {\rm d} t \le c.
\end{align*}
Then the definitions of $w$ and the Sobolev--Slobodeckij seminorm lead to
    \begin{align*}
     & |w|_{W^{\alpha-\epsilon,1}(0,T;\dot{H}^{-1-\mu}(\Omega))}  = \int_0^T\hskip-5pt\int_0^T \frac{\|w(t)-w(s)\|_{\dot{H}^{-1-\mu}(\Omega)}}{|t-s|^{1+\alpha-\epsilon}}\ {\rm d} t{\rm d}s \\
     \leq &\int_0^T\hskip-5pt\int_0^T \frac{|\int_s^t\|w'(\zeta)\|_{\dot{H}^{-1-\mu}(\Omega)}\ {\rm d}\zeta|}{|t-s|^{1+\alpha-\epsilon}}\ {\rm d} t{\rm d}s
      \leq \int_0^T\hskip-5pt\int_0^T \frac{|t^{\alpha-1}-s^{\alpha-1}|}{|t-s|^{1+\alpha-\epsilon}}\ {\rm d} t{\rm d}s \\
      =& c\int_0^1\hskip-5pt\int_0^1 \frac{|\xi^{\alpha-1}-\eta^{\alpha-1}|}{|\xi-\eta|^{1+\alpha-\epsilon}}\ {\rm d} \eta{\rm d}\xi
      = 2c\int_0^1\hskip-5pt\int_0^\xi \frac{\eta^{\alpha-1}-\xi^{\alpha-1}}{(\xi-\eta)^{1+\alpha-\epsilon}}\ {\rm d} \eta{\rm d}\xi  \\
     =& 2c\int_{0}^{1}\xi^{-1+\epsilon}\ {\rm d} \xi \int_{0}^{1} \frac{ t^{\alpha-1}-1}{(1-t)^{1+\alpha-\epsilon}} \ {\rm d} t \leq c.
    \end{align*}
    This completes the proof of the lemma.
\end{proof}

\section{High-order time-stepping scheme and error analysis}\label{sec:err}

Now we develop a corrected BDF2 time-stepping scheme for problem \eqref{eqn:fde}, and derive an $L^2(\Omega)$ error bound.

\subsection{High order time-stepping scheme}
First, we divide the time interval $[0,T]$ into $N$ uniform subintervals, with a time step size $\tau=T/N$ and grid points $t_n=n\tau$, $n=0,1,\dots,N$. To approximate the Caputo fractional derivative $\Dal v $, we employ BDF2 convolution quadrature defined by
\begin{equation*}
	\bar{\partial}_\tau^\alpha (v(t_n)-v_0) :=\tau^{-\alpha}\sum_{j=0}^nb^{(\alpha)}_{j}(v(t_{n-j})-v_0),\quad\mbox{with}\ v_0=v(0),
\end{equation*}
where the weights $b_j^{(\alpha)}$ are generated by the power series expansion
\begin{equation*}
  \delta(\zeta)^\alpha=\sum_{j=0}^\infty b_j^{(\alpha)}\zeta^j,\quad \mbox{with}\quad \delta(\zeta):=\frac12\zeta^2-2\zeta+\frac32.
\end{equation*}

Then we defined an implicit time semi-discrete scheme for problem
\eqref{eqn:fde}: Find a sequence $(u^n)_{n=0}^N\subset H_0^1(\Omega)$ with $u^0=u_0$ such that
\begin{equation}\label{eqn:BDF2-scheme for u}
	\left\{\begin{aligned}
		\bar{\partial}_\tau^\alpha u^1 + A(u_0)u^1 &= f(t_1)+ \tfrac{1}{2}(f(0)-A(u_0)u_0) + [A(u_0)-A(u^1)]u^1, \\
		\bar{\partial}_\tau^\alpha u^n + A(u_0)u^n  &=  f(t_n)+[A(u_0)-A(u^n)]u^n, \quad n=2,3,\dots,N.
	\end{aligned}
	\right.
\end{equation}
The scheme involves correction at the first step using the term $\frac12(f(0)-A(u_0)u_0)$. The construction is inspired by the scheme in \cite{jin2016two,jin2017correction} for convolution quadrature for linear subdiffusion problems. It also resembles that for linear problems with time-dependent coefficients \cite[Section 9.4]{JinZhou:2023book}. In Theorem \ref{thm:main-err-estimate}, we shall prove that the correction allows achieving high-order convergence $O(\tau^{1+\alpha-\epsilon})$ for small $\epsilon>0$.

The scheme \eqref{eqn:BDF2-scheme for u} involves solving a nonlinear elliptic equation at each time step.
Its unique solvability follows from Brouwer's fixed point theorem \cite[p. 237]{Thome2006GalerkinFE}.
Indeed, by multiplying by $\tau^\alpha$ on both sides and treating the terms involving $u^1,\dots, u^{n-1}$ as the known source, we can define a
continuous mapping $M: L^2(\Omega)\mapsto L^2(\Omega)$ for \eqref{eqn:BDF2-scheme for u} with a zero Dirichlet boundary condition, and the
mapping $M$ satisfies for all $n=1,2,\dots,N$
\begin{equation*}
	(Mu^n,u^n):=([b_0^{(\alpha)}+ \tau^\alpha A(u^n)]u^n,u^n )-\tau^\alpha(w,u^n)\geq b_0^{(\alpha)}\|u^n\|^2_{L^2(\Omega)}-\tau^\alpha\|w\|_{L^2(\Omega)}\|u^n\|_{L^2(\Omega)},
\end{equation*}
for some $w\in L^2(\Omega)$. Then for small enough $\tau$, the strict positivity of the map $M$ follows directly. By Brouwer's fixed point theorem, we deduce the existence of the solution sequence $(u^n)_{n=0}^N$ to the scheme \eqref{eqn:BDF2-scheme for u}. Moreover, we can rearrange \eqref{eqn:BDF2-scheme for u} as a nonlinear elliptic problem
\begin{equation*}
	[\tau^{-\alpha} b_0^{(\alpha)}+A(u^1)]u^1 = \tfrac{1}{2}(f(0)-A(u_0)u_0) + f(t_1) +\tau^{-\alpha} b_1^{(\alpha)}u_0.
\end{equation*}
By Assumption \ref{ass:reg-high} and the regularity result in \cite[Theorem 2.5]{Casas:2009},
we have $u^1\in W^{2,4}(\Omega)\cap H^1_0(\Omega)\hookrightarrow W^{1,\infty}(\Omega)$ for small $\tau$. Then by mathematical induction, we
can deduce that the solution sequence $(u^n)_{n=2}^N$ belongs to $ W^{2,4}(\Omega)\cap H^1_0(\Omega)$. This regularity estimate will play a role in the error analysis.

\subsection{Error analysis}

Now we give an error analysis of the time discrete solution $(u^n)_{n=0}^N$ of the corrected scheme \eqref{eqn:BDF2-scheme for u}. The analysis employs the following splitting of problem \eqref{eqn:fde}: find $v\equiv v(t)\in H_0^1(\Omega)$ with $v(0) = u_0$ and $w\equiv w(t)\in H_0^1(\Omega)$ with $w(0) = 0$ such that for a.e. $t\in(0,T]$,
\begin{equation}\label{eqn: prob-splitting}
	\Dal  v + A(u_0)v = f \quad\mbox{and}\quad
	\Dal  w + A(u_0)w =[ A(u_0)-A(u)] u.
\end{equation}
The function $v$ solves a linearized problem, and the function $w$ solves a linearized problem with the (nonlinear) perturbed remainder.
The corrected BDF2 time-stepping schemes for $v$ and $w$ read: find $(v^n)_{n=0}^N\subset H_0^1(\Omega)$ with  $v^0=u_0$ such that \begin{equation}\label{eqn:BDF2-scheme for v}
	\left\{\begin{aligned}
		\bar{\partial}_\tau^\alpha v^1+ A(u_0)v^1   &= f(t_1) + \tfrac{1}{2}(f(0)-A(u_0)u_0),  \\
		\bar{\partial}_\tau^\alpha v^n + A(u_0)v^n  &=  f(t_n), \quad n=2,3,\dots,N,
	\end{aligned}
	\right.
\end{equation}
and find $(w^n)_{n=0}^N\subset H_0^1(\Omega)$ with  $w^0=0$ such that
\begin{equation}\label{eqn:BDF2-scheme for w}
	\bar{\partial}_\tau^\alpha w^n + A(u_0)w^n = [A(u_0)-A(u^n)]u^n, \quad n=1,2,\dots,N.
\end{equation}
Also we define an auxiliary problem: Find $(\overline{w}^n)_{n=1}^N\subset H_0^1(\Omega)$ with $\overline{w}^0=0$ such that
\begin{equation}\label{eqn:BDF2-scheme for bar w}
\bar{\partial}_\tau^\alpha\overline{w}^n + A(u_0)\overline{w}^n  = [A(u_0)-A(u(t_n))]u(t_n), \quad n=1,2,\dots,N.
\end{equation}
For problem \eqref{eqn:BDF2-scheme for v}, under Assumption \ref{assum: reg-plus}, the following  error bound holds \cite[Theorem 3.4]{JinZhou:2023book}
\begin{equation}\label{ineq:err-v}
	\|v(t_n)-v^n\|_{L^2(\Omega)}\leq c\tau^2t_n^{\alpha-2}, \quad  n=1,2,\dots,N.
\end{equation}

In view of the etimate \eqref{ineq:err-v} and the splitting
\begin{align}\label{eqn:split-time}
    u(t_n) - u^n & = (v(t_n)-v^n) + (w(t_n)-w^n)\nonumber\\
     & = (v(t_n)-v^n) + (w(t_n)-\overline{w}^n) + (\overline{w}^n-w^n),
\end{align}
it suffices to bound the error components $w(t_n)-\overline{w}^n$ and $\overline{w}^n-w^n$. To this end, we define the semidiscrete solution operator $E_\tau^{n}(A(u_0))$ by
\begin{equation*}
	E_\tau^{n}(A(u_0)) = \frac{1}{2\pi i}\int_{\Gamma^\tau_{\theta,\sigma}}e^{zt_n}(\delta_\tau(e^{-z\tau})^\alpha+A(u_0))^{-1}\ {\rm d}z,
\end{equation*}
with $\Gamma_{\theta,\sigma}^\tau:=\{z\in\Gamma_{\theta,\sigma}:\sigma\leq |z|\leq \frac{\pi\sin\theta}{\tau}\}$ oriented with an increasing imaginary part. The following smoothing property of $E_\tau^{n}(A(u_0))$ holds \cite[p. 252, (9.38)]{JinZhou:2023book}: for any $\gamma\in[0,1]$,
\begin{equation}\label{ineq:dis-smooth-prop}
	\|A(u_0)^\gamma E_\tau^{n}(A(u_0))\| \leq c(t_n+\tau)^{(1-\gamma)\alpha-1}, \quad n=0,1,\dots, N.
\end{equation}
Next, we recall several useful estimates  \cite[Lemma 3.1]{JinZhou:2023book}.
\begin{lemma}\label{lem:disc-kernel-basic}
Let $\omega_\tau(e^{-z\tau})=\tau^{-1}(1-e^{-z\tau})$ and $\delta_\tau(e^{-z\tau}) = \tau^{-1}\delta(e^{-z\tau})$. Then for any fixed $\theta\in (\frac{\pi}{2},\pi)$, there exists $\theta'\in(\frac{\pi}{2},\pi)$ such that for all $\alpha\in(0,1)$ and $z\in \Gamma^\tau_{\theta,\sigma}$, we have $\omega_\tau(e^{-z\tau})$, $\delta_\tau(e^{-z\tau})\in\Sigma_{\theta'}$, and moreover,
\begin{align}
c_1|z|\leq|\omega_\tau(e^{-z\tau})|\leq c_2|z|,\
|\omega_\tau(e^{-z\tau})-z|\leq c_3\tau|z|^{2}, \ |\omega_\tau(e^{-z\tau})^\alpha-z^\alpha|\leq c_3\tau|z|^{1+\alpha}, \label{eqn:disc-kern-1} \\
c_1|z|\leq|\delta_\tau(e^{-z\tau})|\leq c_2|z|,\ |\delta_\tau(e^{-z\tau})-z|\leq c_3\tau^2|z|^{3}, \ |\delta_\tau(e^{-z\tau})^\alpha\!-\!z^\alpha|\leq c_3\tau^2|z|^{2+\alpha}. \label{eqn:disc-kern-2}
\end{align}
\end{lemma}

The next lemma gives an $L^2(\Omega)$ error bound on $\bar{w}^n$.
\begin{lemma}\label{lem:err-inter for bar w}
	Let $\mu>0$ be small, and Assumption \ref{assum: reg-plus} hold. Let $w\equiv w(t)$ and $(\overline {w}^n)_{n=1}^N$ solve problems \eqref{eqn: prob-splitting} and \eqref{eqn:BDF2-scheme for bar w}, respectively. Then there exists $c>0$, independent of $\tau$, such that
	\begin{equation*}
		\|w(t_n)-\overline{w}^n\|_{L^2(\Omega)}\leq c\tau^{1+\alpha}t_n^{(\frac12-\mu)\alpha-1}, \quad \ n=1,2,\dots,N.
	\end{equation*}
\end{lemma}
\begin{proof}
For the nonlinear remainder $R(t)\equiv R(u;u_0)$ in \eqref{eqn:remain}, we have
\begin{equation*}
	R(t) = \int_0^t R'(s)\ {\rm d}s =  \Big( \int_0^\tau R'(s) \ {\rm d}s - \tau R'(\tau)\Big) + tR'(\tau) + \Big[ \int_\tau^t(t-s)R''(s)\ {\rm d}s\Big] := \sum_{i=1}^3 {\rm I}_i.
\end{equation*}
 By Laplace transform and generating function, the errors $e_i$, $i=1,2,3$, can be represented by
    \begin{equation*}
    	e_i := \frac{1}{2\pi i}\int_{\Gamma_{\theta,\sigma}}e^{zt_n} (z^\alpha+A(u_0))^{-1} \widehat{{\rm I}}_i \ {\rm d}z - \frac{\tau}{2\pi i}\int_{\Gamma^\tau_{\theta,\sigma}}e^{zt_n}(\delta_\tau(e^{-z\tau})^\alpha+A(u_0))^{-1} \widetilde{{\rm I}}_i(e^{-z\tau}) \ {\rm d}z.
    \end{equation*}
where $~\widehat{\cdot}~$ denotes taking Laplace transform and for a sequence $g = (g^j)_{j=0}^\infty $, $~\widetilde{g}~$ denotes its generating function $\widetilde{g}(\zeta):=\sum_{j=0}^\infty g^j\zeta^j$. Below we bound $e_i$ separately. We define the kernel $K(z)$ by
\begin{equation*}
K(z):=z^{-1}(z^\alpha+A(u_0))^{-1} - \omega_\tau(e^{-z\tau})^{-1}(\delta_\tau(e^{-z\tau})^\alpha+A(u_0))^{-1},
\end{equation*}
which upon rearranging the terms can be written as
\begin{align*}
K(z) =& (z^{-1}-\omega_\tau(e^{-z\tau})^{-1})(z^\alpha+A(u_0))^{-1} \\
&+\omega_\tau(e^{-z\tau})^{-1}(\delta_\tau(e^{-z\tau})^\alpha-z^\alpha)(z^\alpha+A(u_0))^{-1}(\delta_\tau(e^{-z\tau})^\alpha+A(u_0))^{-1}.
\end{align*}
Then Lemma \ref{lem:disc-kernel-basic} and the resolvent estimate \eqref{eqn:resolvent} (with $q=2$) imply
\begin{equation*}
	\|A(u_0)^{\frac12} K(z)\|_{\mathcal{B}(L^2)} \leq c\tau|z|^{-\frac{\alpha}{2}} + c\tau^2|z|^{1-\frac{\alpha}{2}}\leq c\tau|z|^{-\frac{\alpha}{2}},\quad \forall z\in\Gamma_{\theta,\sigma}^\tau.
\end{equation*}
Meanwhile, Theorem \ref{thm:reg-rem} yields \begin{equation}\label{eqn:err-I1}
    \|A(u_0)^{-\frac12}{\rm I}_1\|_{L^2(\Omega)}\leq c\tau^\alpha.
\end{equation}
Consequently,
\begin{equation*}
		\Big\|\frac{1}{2\pi i}\int_{\Gamma^\tau_{\theta,\sigma}}e^{zt_n} K(z) {\rm I}_1 \ {\rm d}z\Big\|_{L^2(\Omega)}
		\leq  c\tau^{1+\alpha}\int_{\Gamma^\tau_{\theta,\sigma}}|e^{zt_n}||z|^{-\frac{\alpha}{2}}\ |{\rm d}z| \leq c\tau^{1+\alpha}t_n^{\frac{\alpha}{2}-1}.
	\end{equation*}
Similarly, the resolvent estimate \eqref{eqn:resolvent} (with $q=2$), Theorem \ref{thm:reg-rem} and the estimate \eqref{eqn:err-I1} lead to
\begin{align*}
& \Big\|\frac{1}{2\pi i}\int_{\Gamma_{\theta,\sigma}\backslash\Gamma^\tau_{\theta,\sigma}}e^{zt_n} z^{-1}(z^\alpha+A(u_0))^{-1}{\rm I}_1\ {\rm d}z\Big\|_{L^2(\Omega)} \\ \leq&
c \int_{\Gamma_{\theta,\sigma}\backslash\Gamma^\tau_{\theta,\sigma}}|e^{zt_n}||z|^{-1}\|A(u_0)^{\frac12}(z^\alpha+A(u_0))^{-1}\|_{\mathcal{B}(L^2)} \|A(u_0)^{-\frac12}{\rm I}_1\|_{L^2(\Omega)} \ |{\rm d}z| \\
    \leq&c\tau^{\alpha} \int^\infty_{\frac{\pi\sin\theta}{\tau}}e^{-crt_n}r^{-1-\frac{\alpha}{2}}\ {\rm d}r\leq c\tau^{1+\alpha}t_n^{\frac{\alpha}{2}-1}.
\end{align*}
Thus, $\|e_1\|_{L^2(\Omega)}\leq c\tau^{1+\alpha}t_n^{\frac{\alpha}{2}-1}$.
Next, let $\gamma(\zeta):= \sum_{n=1}^\infty n\zeta^n = (\zeta\frac{\rm d}{{\rm d}\zeta})\frac{1}{1-\zeta}$ and
\begin{align*}
	\widetilde K(z):=& z^{-2}(z^\alpha+A(u_0))^{-1} - \gamma(e^{-z\tau})\tau^2(\delta_\tau(e^{-z\tau})^\alpha+A(u_0))^{-1}\\
 = & (z^{-2}-\gamma(e^{-z\tau})\tau^2) (\delta_\tau(e^{-z\tau})^\alpha+A(u_0))^{-1} \\
   &+  z^{-2}(\delta_\tau(e^{-z\tau})^\alpha-z^\alpha)(z^\alpha+A(u_0))^{-1}(\delta_\tau(e^{-z\tau})^\alpha+A(u_0))^{-1}.
\end{align*}
From the estimate $|\gamma(e^{-z\tau})\tau^2-z^{-2}|\leq c\tau^2$ \cite[p. 60]{JinZhou:2023book}, \eqref{eqn:disc-kern-2} and the resolvent estimate \eqref{eqn:resolvent}, it follows that
	\begin{equation*}
		\|A(u_0)^{\frac12}\widetilde K(z)\| \leq c\tau^2|z|^{-\frac{\alpha}{2}}, \quad\forall z\in\Gamma_{\theta,\sigma}^\tau.
\end{equation*}
On the contour $\Gamma_{\theta,\sigma}\backslash\Gamma^\tau_{\theta,\sigma}$, the resolvent estimate \eqref{eqn:resolvent} and Theorem \ref{thm:reg-rem} yield
\begin{align*}
	&\Big\|\frac{1}{2\pi i}\int_{\Gamma_{\theta,\sigma}\backslash\Gamma^\tau_{\theta,\sigma}}e^{zt_n} z^{-2}(z^\alpha+A(u_0))^{-1} R'(\tau)\ {\rm d}z\Big\|_{L^2(\Omega)} \\
	\leq&c \int_{\Gamma_{\theta,\sigma}\backslash\Gamma^\tau_{\theta,\sigma}}|e^{zt_n}|| z|^{-2}\|A(u_0)^{\frac12}(z^\alpha+A(u_0))^{-1}\|_{\mathcal{B}(L^2)} \|A(u_0)^{-\frac12}R'(\tau)\|_{L^2(\Omega)} \ |{\rm d}z| \\
			\leq&c\tau^{\alpha-1} \int^\infty_{\frac{\pi\sin\theta}{\tau}}e^{-crt_n}r^{-2-\frac{\alpha}{2}}\ {\rm d}r\leq c\tau^{1+\alpha}t_n^{\frac{\alpha}{2}-1}.
\end{align*}
Repeating the argument for the error $e_1$ gives
\begin{align*}
	&\Big\|\frac{1}{2\pi i}\int_{\Gamma^\tau_{\theta,\sigma}}e^{zt_n}\widetilde K(z) R'(\tau)\ {\rm d}z \Big\|_{L^2(\Omega)} \\
 \leq& c\int_{\Gamma^\tau_{\theta,\sigma}}|e^{zt_n}|\|A(u_0)^{\frac12}\widetilde K(z)\|_{\mathcal{B}(L^2)} \|A(u_0)^{-\frac12}R'(\tau)\|_{L^2(\Omega)} \ |{\rm d}z| \leq c\tau^{1+\alpha}t_n^{\frac{\alpha}{2}-1}.
\end{align*}
Thus we obtain
$\|e_2\|_{L^2(\Omega)}\leq c\tau^{1+\alpha}t_n^{\frac{\alpha}{2}-1}$.
Next, let $E_\tau(t;A(u_0))=\tau\sum_{n=0}^\infty E_\tau^n(A(u_0))\delta_{t_n}(t)$, with $\delta_{t_n}(t)$ being the Dirac-delta function. By the solution representation, we have
\begin{align*}
	e_3 &= \int_\tau^{t_n}\big[E(t_n-s;A(u_0)) - E_\tau(t_n-s;A(u_0))\big]\int_\tau^s(s-\xi)R''(\xi)\ {\rm d}\xi\ {\rm d}s.
\end{align*}
Then for small $\mu>0$, Fubini's theorem implies
\begin{align*}
    \|e_3\|_{L^2(\Omega)} \leq&\int_\tau^{t_n}\|A(u_0)^{\frac12+\frac{\mu}{2}}(E(t_n-s;A(u_0)) - E_\tau(t_n-s;A(u_0))\|_{\mathcal{B}(L^2)}\\
    & \times \int_\tau^s(s-\xi)\|A(u_0)^{-\frac12-\frac{\mu}{2}}R''(\xi)\|_{L^2(\Omega)}\ {\rm d}\xi\ {\rm d}s\\
    =&\int_\tau^{t_n} \|A(u_0)^{-\frac12-\frac{\mu}{2}}R''(\xi)\|_{L^2\II}{\rm II}(t_n-\xi)\ {\rm d}\xi,
\end{align*}
with the factor ${\rm II}(t_n-\xi)$ given by
$${\rm II}(t_n-\xi)=\int_0^{t_n-\xi}\|A(u_0)^{\frac12+\frac{\mu}{2}}[E(t_n-\xi-s;A(u_0)) - E_\tau(t_n-\xi-s;A(u_0))]\|_{\mathcal{B}(L^2)}s\ {\rm d}s.$$
The argument of  \cite[Theorem 3.4]{JinZhou:2023book} leads to
${\rm II}(t_n-\xi) \leq c(t_{n+1}-\xi)^{(\frac12-\frac{\mu}{2})\alpha-1}\tau^2.$
Consequently, Theorem \ref{thm:reg-rem} yields
\begin{align*}
    \|e_3\|_{L^2(\Omega)} \leq c\tau^2 \int_\tau^{t_n} \xi^{\alpha-2}(t_{n+1}-\xi)^{(\frac12-\frac{\mu}{2})\alpha-1}\ {\rm d}\xi \leq c\tau^{1+\alpha}t_n^{(\frac12-\frac{\mu}{2})\alpha-1} + c\tau^2t_n^{\alpha-2}.
\end{align*}
The preceding estimates together show the desired assertion.
\end{proof}

Next we give an $L^2\II$ bound on the error $\overline{w}^n-w^n$.
\begin{lemma}\label{lem:err-inter for bar w and w}
Let the conditions in Lemma \ref{lem:err-inter for bar w} hold. Let $(w^n)_{n=0}^N$ and $(\overline{w}^n)_{n=0}^N$ solve problems \eqref{eqn:BDF2-scheme for w} and \eqref{eqn:BDF2-scheme for bar w}, respectively. Then with $\ell_\tau:=|\log\tau|$, the following estimate holds for $n = 1,\dots,N$:
\begin{equation*}
	\|\overline{w}^n-w^n\|_{L^2(\Omega)} \leq c\tau\sum_{j=1}^n  t_{n-j+1}^{\frac{\alpha}{2}-1}\|u(t_j)-u^j\|_{L^2(\Omega)} + c\ell_\tau\tau^{1+\alpha}t_n^{(\frac12-\mu)\alpha-1}.
\end{equation*}
\end{lemma}
\begin{proof}
Let $e^n:=\overline{w}^n-w^n$ with $e^0=0$. By \eqref{eqn:BDF2-scheme for w} and \eqref{eqn:BDF2-scheme for bar w} and the splitting \eqref{eqn:split-time}, the function $e^n$ satisfies for $n = 1,2,\dots,N$:
\begin{align*}
\bar{\partial}_\tau^\alpha e^n + A(u(t_m))e^n
=& [A(u_0)-A(u(t_n))]u(t_n) - [A(u_0)-A(u^n)]u^n \\
 &- [A(u_0)-A(u(t_m))](\overline{w}^n - w^n)\\
=& [A(u(t_m))-A(u(t_n))]u(t_n) - [A(u(t_m))-A(u^n)]u^n\\
 &+  [A(u_0)-A(u(t_m))](v(t_n)-v^n + w(t_n)-\overline{w}^n) \\
=& [A(u(t_m))-A(u(t_n))](u(t_n)-u^n) - [A(u(t_n))-A(u^n)]u^n \\
 &+  [A(u_0)-A(u(t_m))](v(t_n)-v^n + w(t_n)-\overline{w}^n)\\
:=& {\rm I}_1^n + {\rm I}_2^n + {\rm I}_3^n.
\end{align*}
By the discrete Laplace transform, the three error components $e_i$ are given by
\begin{equation*}
	e^m_i:= \tau\sum_{j=0}^mE_\tau^{m-j}(A(u(t_m))) {\rm I}_i^j,\quad i=1,2,3.
\end{equation*}
For the term $e_1^m$, the estimates \eqref{ineq:sharp-holder-perturb} and \eqref{ineq:dis-smooth-prop} yield
\begin{align*}
&\|e_1^m\|_{L^2(\Omega)} \\
 \leq &\tau\sum_{j=0}^m\|A(u(t_m))E_\tau^{m-j}(A(u(t_m)))\|_{\mathcal{B}(L^2)} \|[I-A(u(t_m))^{-1}A(u(t_j))](u(t_j)- u^j)\|_{L^2(\Omega)} \\
  \leq& c\tau\sum_{j=1}^mt_{m-j+1}^{-1}t_{m-j}^{\frac{\alpha}{2}}\|u(t_j)-u^j\|_{L^2(\Omega)} \leq c\tau\sum_{j=1}^m t_{m-j+1}^{\frac{\alpha}{2}-1}\|u(t_j)-u^j\|_{L^2(\Omega)}.
\end{align*}
H\"{o}lder's inequality, the regularity $(u^n)_{n=0}^N\subset W^{1,\infty}(\Omega)$ and Lemma \ref{lem:Riesz} lead to
\begin{align*}
    &\|A(u(t_m))^{-\frac12}[A(u(t_j))-A(u^j)]u^j\|_{L^2(\Omega)}\\
    =& \sup_{\|\varphi\|_{L^2(\Omega)=1}} ([a(u(t_j))-a(u^j)]\nabla u^j,\nabla A(u(t_m))^{-\frac12}\varphi) \\
    \leq&c\|u(t_j)-u^j\|_{L^2(\Omega)}\|\nabla u^j\|_{L^\infty(\Omega)}\|\nabla A(u(t_m))^{-\frac12}\|_{\mathcal{B}(L^2)}\\
    \leq& c\|u(t_j)-u^j\|_{L^2(\Omega)}.
\end{align*}	
This and the estimate \eqref{ineq:dis-smooth-prop} imply
\begin{align*}
  &\|e_2^m\|_{L^2(\Omega)}  \\
  \leq& \tau\sum_{j=0}^m\|A(u(t_m))^\frac{1}{2}E_\tau^{m-j}(A(u(t_m)))\|_{\mathcal{B}(L^2)} \|A(u(t_m))^{-\frac12}[A(u(t_j))-A(u^j)]u^j\|_{L^2(\Omega)} \\
  \leq&c\tau\sum_{j=1}^mt_{m-j+1}^{\frac{\alpha}{2}-1}\|u(t_j)-u^j\|_{L^2(\Omega)}.
\end{align*}
Using the estimates \eqref{ineq:err-v} and \eqref{ineq:dis-smooth-prop}, and Lemma \ref{lem:err-inter for bar w}, we derive
\begin{align*}
    \|e_3^m\|_{L^2(\Omega)} & \leq \tau\sum_{j=1}^m\|A(u(t_m))E_\tau^{m-j}(A(u(t_m)))\|_{\mathcal{B}(L^2)} \\
      &\quad\times \|[A(u(t_m))^{-1}A(u_0)-I](v(t_j)-v^j + w(t_j)-\overline{w}^j)\|_{L^2(\Omega)} \\
    & \leq ct_m^{\frac{\alpha}{2}}\tau\sum_{j=1}^mt_{m-j+1}^{-1}\big[\|v(t_j)-v^j \|_{L^2(\Omega)} + \|w(t_j)-\overline{w}^j \|_{L^2(\Omega)}\big]\\
    & \leq c_T(\tau^3\sum_{j=1}^mt_{m-j+1}^{-1}t_j^{\alpha-2}+ \tau^{2+\alpha} \sum_{j=1}^mt_{m-j+1}^{-1}t_j^{(\frac12-\mu)\alpha-1}) \leq c\ell_\tau\tau^{1+\alpha}t_m^{(\frac12-\mu)\alpha-1}.
\end{align*}
Since the assertion holds for any $1\leq m\leq n \leq N$, the desired estimate follows directly.
\end{proof}

Below we need a useful Gronwall’s inequality, which generalizes the standard variants in \cite[Lemma 7.1]{ElliottLarsson:1992} and \cite[Lemma 9.9]{JinZhou:2023book}.
\begin{lemma}\label{lem:gronwalls}
Let $ 0\le \varphi^n \leq R$ for $0 \leq t_n \leq T$. With $\ell_\tau:=|\log\tau|$, if
$$
\varphi^n \leq a t_n^{-1} \ell_\tau^p+b \tau \sum_{j=1}^n t_{n-j+1}^{\beta-1}\varphi^j, \quad 0<t_n \leq T,
$$
for some $a, b \geq 0$, $\beta\in(0,1]$ and $p>0$, then there is $c=c(b, \beta, T, R)$ and $k=\lfloor 1/\beta \rfloor$ such that
$$
\varphi^n \leq  cat_n^{-1} \ell_\tau^{p+k}, \quad 0<t_n \leq T.
$$
\end{lemma}

\begin{proof}
Using the estimate
$\tau  \sum_{j=1}^n t_{n-j+1}^{-\beta} t_j^{-1} \le c \tau^{-\beta} \ell_\tau$
(see e.g., \cite[Lemma 3.11]{JinZhou:2023book}) and iterating the given inequality $k$ times yield
\begin{equation*}
\varphi^n \leq c a t_n^{-1} \ell_\tau^{p+k} + c \tau \sum_{j=1}^n t_{n-j+1}^{(k+1)\beta-1}\varphi^j
\le c a t_n^{-1} \ell_\tau^{p+k}
+  c \tau \sum_{j=1}^n \varphi^j.
\end{equation*}
Now the assertion follows from the discrete Gronwall's inequality with log factors \cite[Lemma 9.9]{JinZhou:2023book}.
\end{proof}

Last we give an $L^2(\Omega)$ error bound on the scheme \eqref{eqn:BDF2-scheme for u}. The  rate $O(\tau^{1+\alpha}|\log\tau|^k)$ is consistent with the regularity of the remainder $R(u;u_0)$ in Theorem \ref{thm:reg-rem}. This rate is slower than the rate $O(\tau^2)$ for the scheme applied to standard linear subdiffusion \cite{jin2016two}, due to the limited temporal regularity of the solution $u$.

\begin{theorem}\label{thm:main-err-estimate}
Let the conditions in Lemmas \ref{lem:err-inter for bar w} and \ref{lem:err-inter for bar w and w} hold, and $ u(t)$ and $(u^n)_{n=0}^N$ be the solutions to problems \eqref{eqn:fde} and \eqref{eqn:BDF2-scheme for u}, respectively. Then with $k=\lceil 2/\alpha \rceil$, there holds:
\begin{equation*}
\|u(t_n)-u^n\|_{L^2(\Omega)}\leq c \tau^{1+\alpha} \ell_\tau^{k} t_n^{-1}, \quad n=1,2,\dots,N.
\end{equation*}
\end{theorem}
\begin{proof}
The identities \eqref{eqn: prob-splitting} and \eqref{eqn:split-time},  the estimate \eqref{ineq:err-v}, and Lemmas \ref{lem:err-inter for bar w} and \ref{lem:err-inter for bar w and w} yield
\begin{align*}
    \|u(t_n)-u^n\|_{L^2(\Omega)} &\leq  c\big( \tau^2t_n^{\alpha-2} + \ell_\tau\tau^{1+\alpha}t_n^{(\frac12-\mu)\alpha-1}\big) + c\tau\sum_{j=1}^n t_{n-j+1}^{\frac{\alpha}{2}-1}\|u(t_j)-u^j\|_{L^2(\Omega)}\\
    &\leq c \ell_\tau\tau^{1+\alpha}t_n^{-1}  + c\tau\sum_{j=1}^n t_{n-j+1}^{\frac{\alpha}{2}-1}\|u(t_j)-u^j\|_{L^2(\Omega)}.
\end{align*}
The desired assertion follows from the discrete Gronwall's inequality in Lemma \ref{lem:gronwalls}.
\end{proof}

\begin{remark}\label{rem:fast}
Convolution quadrature (CQ) admits fast and memory-efficient implementations. For pioneering works on fast and oblivious quadrature, see \cite{LubichSchadle:2002, SchadleLopezLubich:2006}. Recent developments of fast CQ algorithms focus on linear time-fractional evolution equation \cite{BanjaiLopez:2019, Fischer:2019}, and for adaptations to quasilinear subdiffusion problems, see \cite[Section 4]{lopez2023convolution}.
\end{remark}

\section{Numerical experiments and discussions}\label{sec:numers}

Now we present numerical experiments to complement the error analysis in Section \ref{sec:err}. We measure the accuracy of the time semi-discrete approximation $(u^n)_{n=1}^N\subset H_0^1(\Omega)$ in the $L^2(\Omega)$ error $e^n:=\|u(t_n)-u^n\|_{L^2(\Omega)}$. Throughout, we fix the terminal time $T=0.1$ and space domain $\Omega=B(0,\frac14)$, a ball centered at the origin of radius $\frac14$. Since the exact solution is unavailable for all the examples, we obtain a reference solution on a finer grid with $\tau=7.81\text{e-5}$. For the space approximation, we employ the conforming piecewise linear Galerkin FEM. To solve the nonlinear elliptic equation at each time step, we employ Newton's method. Numerically, the convergence of Newton iteration is achieved within a few iterations, and hence the scheme is computationally efficient. All the experiments are conducted using the public FEM package FreeFEM++ \cite{hecht2012new} on a personal laptop.

\begin{example}\label{exam:2d}
Fix $\alpha\in\{0.15, 0.30,0.45,0.60,0.75, 0.90\}$, and consider the following settings for problem data:
\begin{enumerate}
		\item[{\rm(a)}] $u_0=1-16x^2-16y^2$, $f= x^2+y^2$ and $a(u) = 1 + \frac15\sin(u)$;
		\item[{\rm(b)}] $u_0=\min(1-16x^2-16y^2,\frac34)$, $f= x^2+y^2$ and $a(u) = 1 + \frac15\sin(u)$;
		\item[{\rm(c)}] $u_0=1-16x^2-16y^2$, $f= tx(1-y)$ and $a(u) = 1 + \frac15\sin(u)$;
		\item[{\rm(d)}] $u_0=\min(1-16x^2-16y^2,\frac34)$, $f= tx(1-y)$ and $a(u) = 1 + \frac15\sin(u)$.
	\end{enumerate}
\end{example}

This example covers various scenarios over a planar circular domain, including both smooth and non-smooth problem data, and rotationally and non-rotationally symmetric cases.  The numerical results are presented in Table \ref{tab:err-2d}. For small fractional order $\alpha\in(0,\frac12)$, the convergence behavior of the error $e^n$ agrees well with the theoretical prediction: the error $e^n$ decays to zero at a rate $\mathcal{O}(\tau^{1+\alpha})$ as $\tau\to 0^+$. The empirical convergence rate is slightly higher than the theoretical prediction for $\alpha\in(\frac12,1)$. These observations hold for all four cases, and the empirical rates are largely comparable irrespective of the regularity of the problem data. In sum, these empirical observations confirm the sharpness of the error estimate in Theorem \ref{thm:main-err-estimate}.

\begin{table}[htpb!]
	\centering\setlength{\tabcolsep}{3pt}
	\caption{Numerical results for the four cases of Example \ref{exam:2d}, with $h=$2.42e-2, and the error $e^n$ at $t=T$.\label{tab:err-2d}}
    \begin{threeparttable}
		\subfigure[]{
			\begin{tabular}{c|cccccc}
			\toprule
			$\alpha\backslash \tau$ & 5.00e-3 & 2.50e-3 & 1.25e-3 & 6.25e-4 & 3.13e-4 & rate \\
			\midrule
			$0.15$ & 2.22e-6 & 8.87e-7 & 3.80e-7 & 1.65e-7 & 6.84e-8 & 1.25\\
			\midrule
			$0.30$ & 5.58e-6 & 2.10e-6 & 8.54e-7 & 3.56e-7 & 1.41e-7 & 1.32\\
			\midrule
			$0.45$ & 9.65e-6 & 3.27e-6 & 1.19e-6 & 4.39e-7 & 1.53e-7 & 1.49\\
			\midrule
			$0.60$ & 1.35e-5 & 3.90e-6 & 1.16e-6 & 3.45e-7 & 9.65e-8 & 1.78\\
			\midrule
			$0.75$ & 1.66e-5 & 4.23e-6 & 1.08e-6 & 2.76e-7 & 6.74e-8 & 1.98\\
			\midrule
			$0.90$ & 1.85e-5 & 4.48e-6 & 1.09e-6 & 2.68e-7 & 6.35e-8 & 2.04\\
			\bottomrule
		\end{tabular}}
	\subfigure[]{
		\begin{tabular}{c|cccccc}
		\toprule
		$\alpha\backslash \tau$ & 5.00e-3 & 2.50e-3 & 1.25e-3 & 6.25e-4 & 3.13e-4 & rate \\
		\midrule
		$0.15$ & 1.85e-6 & 7.27e-7 & 3.08e-7 & 1.33e-7 & 5.49e-8 & 1.27 \\
		\midrule
		$0.30$ & 4.70e-6 & 1.73e-6 & 6.97e-7 & 2.88e-7 & 1.13e-7 & 1.34 \\
		\midrule
		$0.45$ & 8.24e-6 & 2.75e-6 & 9.88e-7 & 3.62e-7 & 1.26e-7 & 1.50 \\
		\midrule
		$0.60$ & 1.17e-5 & 3.36e-6 & 1.00e-6 & 2.96e-7 & 8.33e-8 & 1.78 \\
		\midrule
		$0.75$ & 1.48e-5 & 3.76e-6 & 9.67e-7 & 2.47e-7 & 6.05e-8 & 1.98 \\
		\midrule
		$0.90$ & 1.67e-5 & 4.04e-6 & 9.88e-7 & 2.42e-7 & 5.74e-8 & 2.04\\
		\bottomrule
	\end{tabular}}
    \subfigure[]{
    	\begin{tabular}{c|cccccc}
    	\toprule
    	$\alpha\backslash \tau$ & 5.00e-3 & 2.50e-3 & 1.25e-3 & 6.25e-4 & 3.13e-4 & rate \\
    	\midrule
    	$0.15$ & 2.38e-6 & 9.83e-7 & 4.32e-7 & 1.91e-7 & 7.98e-8 & 1.22 \\
    	\midrule
    	$0.30$ & 5.65e-6 & 2.14e-6 & 8.79e-7 & 3.69e-7 & 1.47e-7 & 1.31 \\
    	\midrule
    	$0.45$ & 9.69e-6 & 3.30e-6 & 1.20e-6 & 4.49e-7 & 1.58e-7 & 1.48 \\
    	\midrule
    	$0.60$ & 1.34e-5 & 3.93e-6 & 1.18e-6 & 3.58e-7 & 1.04e-7 & 1.75 \\
    	\midrule
    	$0.75$ & 1.66e-5 & 4.25e-6 & 1.10e-6 & 2.92e-7 & 7.89e-8 & 1.93 \\
    	\midrule
    	$0.90$ & 1.86e-5 & 4.50e-6 & 1.11e-6 & 2.84e-7 & 7.56e-8 & 1.98 \\
    	\bottomrule
    \end{tabular}}
    \subfigure[]{
    	\begin{tabular}{c|cccccc}
    		\toprule
    		$\alpha\backslash \tau$ & 5.00e-3 & 2.50e-3 & 1.25e-3 & 6.25e-4 & 3.13e-4 & rate \\
    		\midrule
    		$0.15$ & 2.04e-6 & 8.41e-7 & 3.70e-7 & 1.64e-7 & 6.85e-8 & 1.22\\
    		\midrule
    		$0.30$ & 4.78e-6 & 1.78e-6 & 7.27e-7 & 3.04e-7 & 1.21e-7 & 1.32\\
    		\midrule
    		$0.45$ & 8.29e-6 & 2.78e-6 & 1.00e-6 & 3.74e-7 & 1.32e-7 & 1.49\\
    		\midrule
    		$0.60$ & 1.17e-5 & 3.39e-6 & 1.02e-6 & 3.11e-7 & 9.29e-8 & 1.74\\
    		\midrule
    		$0.75$ & 1.48e-5 & 3.78e-6 & 9.89e-7 & 2.65e-7 & 7.30e-8 & 1.91\\
    		\midrule
    		$0.90$ & 1.67e-5 & 4.06e-6 & 1.01e-6 & 2.60e-7 & 7.06e-8 & 1.97\\
    		\bottomrule
    \end{tabular}}
	\end{threeparttable}
\end{table}

\begin{example}\label{exam:1d}
Fix $\alpha\in\{0.30,0.60,0.90\}$, and take $T=10.0$, $\Omega = (0,1)$, $u_0=1-x^2$, $f= 4x^2$ and $a(u) = 1 + \frac15\sin(u)$. The boundary conditions are given by $\frac{{\rm d}}{{\rm d}x}u(0,t) = u(1,t)=0$ for all $t\in[0,T]$.
\end{example}

\begin{figure}[htbp!]
\centering
\setlength{\tabcolsep}{0pt}
\begin{tabular}{ccc}
\includegraphics[width=0.33\textwidth,trim={3.5cm 9cm 4cm 9cm}, clip]{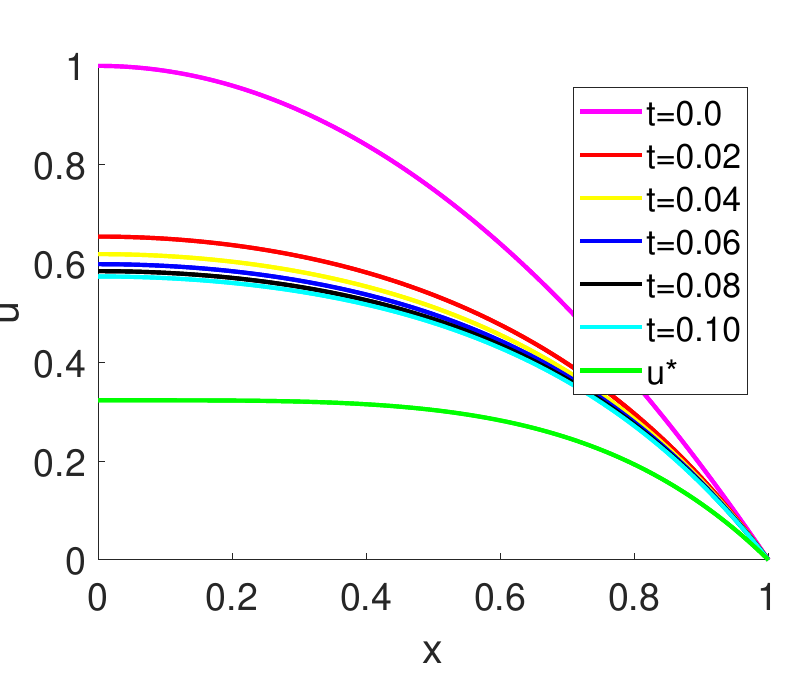} & \includegraphics[width=0.33\textwidth,trim={3.5cm 9cm 4cm 9cm}, clip]{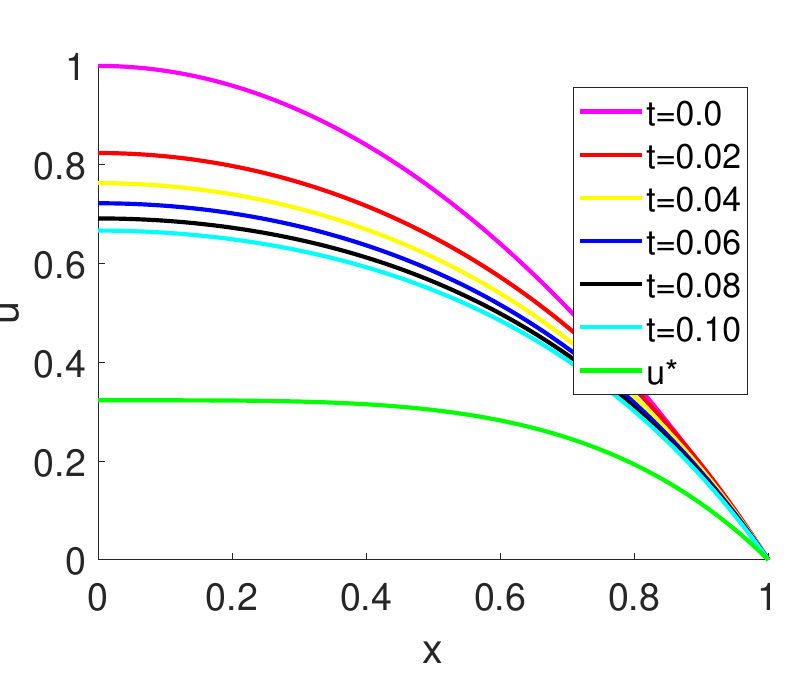} & \includegraphics[width=0.33\textwidth,trim={3.5cm 9cm 4cm 9cm}, clip]{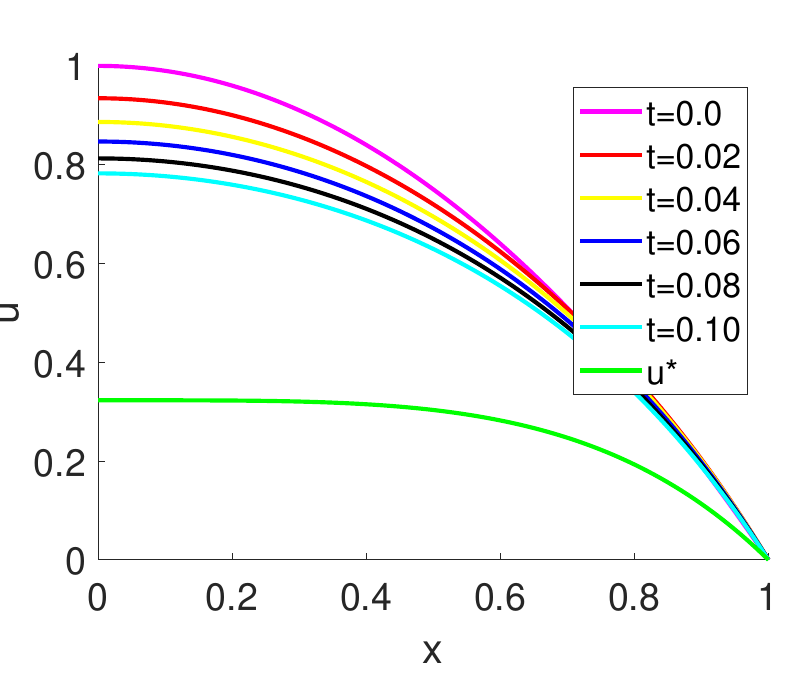}\\
\includegraphics[width=0.33\textwidth,trim={3.5cm 9cm 4cm 9cm}, clip]{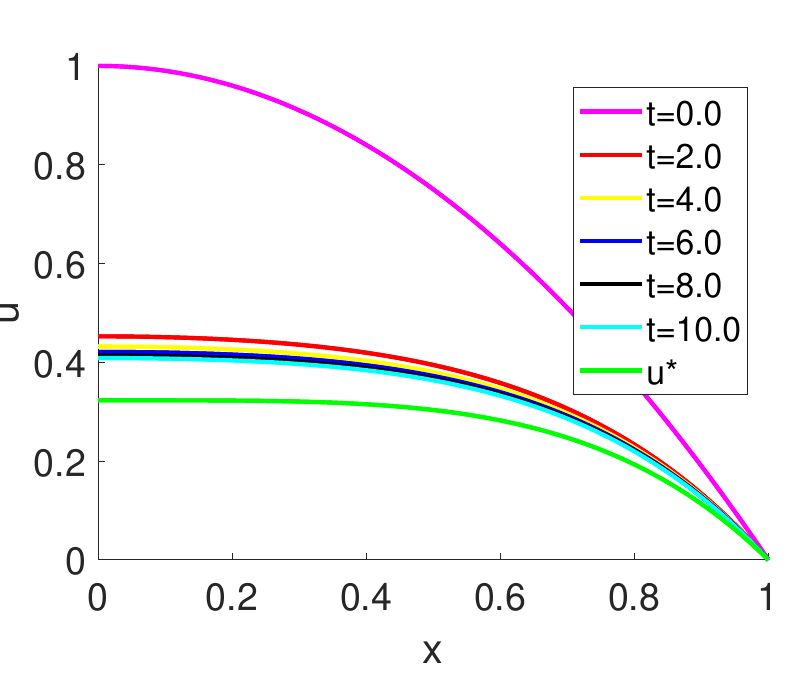} & \includegraphics[width=0.33\textwidth,trim={3.5cm 9cm 4cm 9cm}, clip]{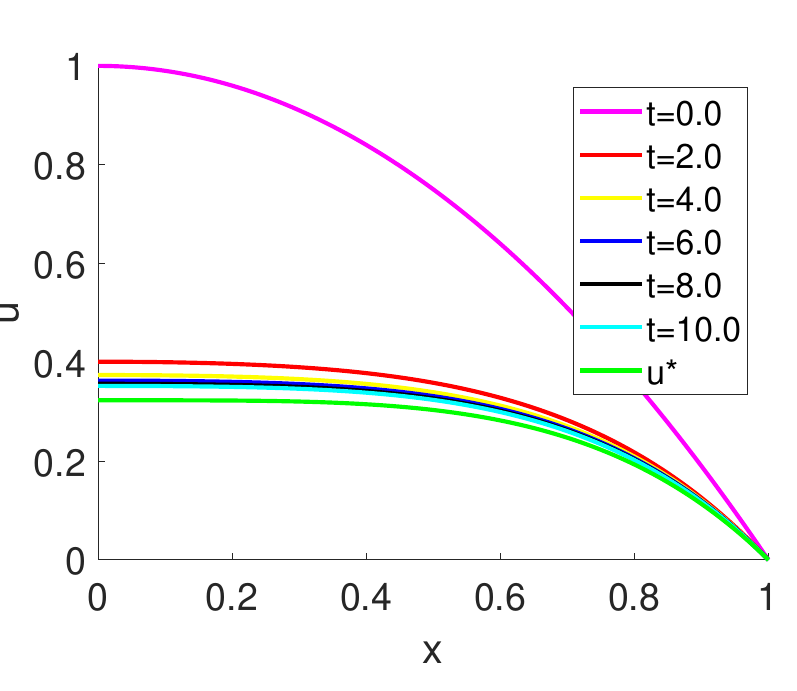} & \includegraphics[width=0.33\textwidth,trim={3.5cm 9cm 4cm 9cm}, clip]{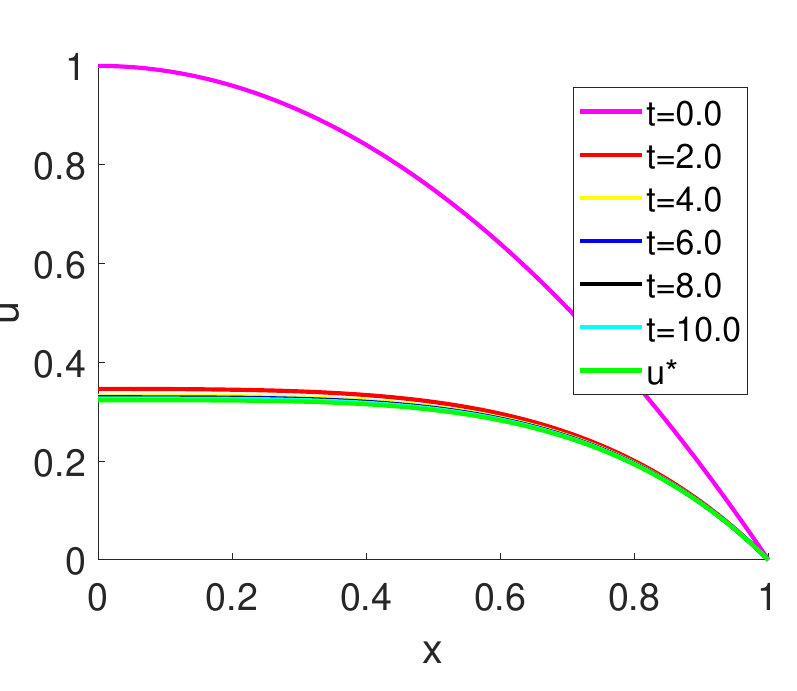}\\

(a) $\alpha = 0.30$ & (b) $\alpha = 0.60$ & (c) $\alpha = 0.90$
\end{tabular}
\caption{Initial condition, solutions at different $t$ and stationary solution $u^*$.}\label{fig:frac1d}
\end{figure}

Example \ref{exam:1d} presents a polar coordinate version of Example \ref{exam:2d}, specifically focusing on the rotationally symmetric case without an angular dependence. The solution profiles in Fig. \ref{fig:frac1d} clearly show that the solution $u (t)$ decays to the steady-state solution $u^*$ as the time $t$ evolves. The smaller is the fractional order $\alpha$, the quicker is the decay of the solution at the initial time. In contrast, when the evolution time $t$ is relatively large, a larger $\alpha$ leads to a faster decay towards the steady-state solution.

\section{Concluding remarks}
In this work we have developed several new pointwise-in-time regularity estimates in Sobolev spaces as well as a high-order time stepping scheme (based on convolution quadrature generated by the second-order backward differentiation formula and suitable correction at the first step) for a quasilinear subdiffusion model, and established nearly optimal error bounds on the scheme, which are confirmed by several numerical experiments with smooth and non-smooth problem data. Note that within the context of subdiffusion, there are alternative models, e.g., $\partial_t u - \nabla\cdot(a(u)\nabla \partial_t^{1-\alpha} u)=f$ (i.e., the elliptic operator acting on the fractional derivative) \cite{Fritz:2024,FritzSuli:2024}. It is of much interest to extend the analysis to these alternatives. Additionally, developing higher-order schemes for the quasilinear subdiffusion model, especially for the case of nonsmooth data, remains crucial. This is technically  very challenging even for the semilinear problem \cite{LiMa:SINUM2022} and problems with time-dependent diffusion coefficients \cite{JinLiZhou:NM2020,Sammon:1983}.

\bibliographystyle{abbrv}
\bibliography{reference}
\end{document}